\documentclass[10pt,reqno]{amsart}
\usepackage{amsmath}
\usepackage{amsthm}
\usepackage{graphicx}
\usepackage{amsfonts}
\usepackage{amssymb}
\usepackage{enumerate}
\usepackage[margin=1.1in]{geometry}
\usepackage{comment}

\usepackage[english]{babel}
\usepackage[latin1]{inputenc}
\usepackage[numbers]{natbib}

\usepackage[normalem]{ulem}

\usepackage[unicode]{hyperref}

\DeclareOldFontCommand{\bf}{\normalfont\bfseries}{\mathbf}
\DeclareOldFontCommand{\cal}{\normalfont\bfseries}{\mathcal}

\newtheorem{theorem}{Theorem}[section]
\newtheorem{lemma}[theorem]{Lemma}
\newtheorem{proposition}[theorem]{Proposition}
\newtheorem{corollary}[theorem]{Corollary}

\theoremstyle{definition}

\newtheorem{remark}[theorem]{Remark}

\setcitestyle{numbers,open={[},close={]}}

\hypersetup{colorlinks, linkcolor={red}, citecolor={green}, urlcolor={blue}}

\newcommand{\diff}{\,\mathrm{d}}

\newcommand\cC{\mathcal C}

\newcommand\cF{\mathcal F}

\newcommand\cL{\mathcal L}

\newcommand\cP{\mathcal P}
\newcommand\cU{\mathcal U}
\newcommand\cV{\mathcal V}
\newcommand\cW{\mathcal W}

\newcommand\EE{\mathbb E}
\newcommand\FF{\mathbb F}

\newcommand\PP{\mathbb P}

\def \E{\mathbb{E}}
\def \F{\mathbb{F}}
\def \G{\mathbb{G}}
\def \H{\mathbb{H}}

\def \L{\mathbb{L}}

\def \N{\mathbb{N}}

\def \P{\mathbb{P}}
\def \Q{\mathbb{Q}}
\def \R{\mathbb{R}}
\def \S{\mathbb{S}}

\def \W{\mathbb{W}}
\def \X{\mathbb{X}}

\def\Gc{\mathcal{G}}

\newcommand{\x}{\mathbf{x}}

\newcommand{\xdim}{m}
\newcommand{\bmdim}{d}

\DeclareMathOperator*{\esssup}{ess\,sup}

\pdfoutput=1
\usepackage{color}

\begin{document}

%\title{Vanishing viscosity reloaded}
%3)\title{On a non-exponential schilder theorem for McKean-Vlasov diffusions and applications}
%2)\title{Zero-noise limits and control of McKean-Vlasov dynamics}
 %1) \title{A probabilistic approach to the zero-noise limit for PDEs on the Wasserstein space}
 \title[Vanishing viscosity on the Wasserstein space]{A probabilistic approach to vanishing viscosity for PDEs on the Wasserstein space}

\author{Ludovic Tangpi}\footnote{Princeton University, ORFE, ludovic.tangpi@princeton.edu.
Financial support from the NSF grant DMS-2005832 is gratefully acknowledged.}

\date{\today}

\maketitle

\begin{abstract}
In this work we prove an analogue, for partial differential equations on the space of probability measures, of the classical vanishing viscosity result known for equations on the Euclidean space. Our result allows in particular to show that the value function arising in various problems of classical mechanics and games can be obtained as the limiting case of second order PDEs. 
The method of proof builds on stochastic analysis arguments and allows for instance to prove a Freindlin-Wentzell large deviation theorem for McKean-Vlasov equations.\\
\end{abstract}

%\tableofcontents

\section{Introduction}
Consider the Hamilton-Jacobi equation
\begin{equation}
\label{eq:HJ}
	\begin{cases}
		\partial_tv(t,x) + f(t, \partial_x v(t,x)\sigma) = 0\quad \text{in } [0,1)\times\R^\xdim\\
		v(1,x) = F(x)\quad \text{for } x \in \R^\xdim
	\end{cases}
\end{equation}
for two given functions $f:[0,1]\times \R^d\to \R$ and $F:\R^\xdim \to \R$ and a constant non-degenerate matrix $\sigma \in \R^{\xdim\times d}$.
Let the function $g$ be the convex conjugate of $f$, i.e.
\begin{equation*}
	g(t,q) := \sup_{z\in \R^d}\big( q\cdot z - f(t,z) \big).
\end{equation*}
It is well-known that the Hamilton-Jacobi equation \eqref{eq:HJ} characterizes the calculus of variations problem
\begin{equation*}
	v(s,x) = \sup_{\varphi}\bigg(F(\Phi(1)) - \int_s^1g(u,\varphi(u))\diff u\bigg)
\end{equation*}
with the supremum taken over bounded measurable maps $\varphi:[0,1]\to \R^d$ and $\Phi$ satisfies $\diff\Phi(t) = \sigma\varphi(t)\diff t $ and $\Phi(s) =x$.
In general, the function $v$ is not a classical solution of \eqref{eq:HJ} everywhere on $[0,1]\times \R^\xdim$ as it is not everywhere differentiable \cite{Fleming_JMA,Fleming-JDE}.
In fact, \eqref{eq:HJ} is typically understood in the viscosity sense to guarantee that $v$ is the unique solution.
Such viscosity solutions were first obtained by the so-called vanishing viscosity method.
That is, as limit as $n$ goes to infinity of solutions of the second order equation
\begin{equation}
	\begin{cases}
		\partial_tv_n(t,x) +\frac{1}{2n}\partial_{xx}v_n(t,x)+ f(t, \partial_x v_n(t,x)\sigma) = 0\quad \text{in } [0,1]\times\R^\xdim\\
		v_n(t,x) = F(x)\quad \text{for } x \in \R^\xdim.
	\end{cases}
\end{equation}
We refer for instance to the monograph of \citet[Chapters 1 $\&$ 2]{Flem-Soner-second} for a great overview on the topic.
Moreover this approach is also used to speed up numerical approximations of the value function $v$ \cite{Waagan,Cap-Leoni}.
The convergence $v_n\to v$ can be guaranteed under mild growth conditions on $f$ and extends far beyond the above setting.
It has been successfully applied to several areas, perhaps most strikingly in the study of risk sensitive optimal control, see e.g. \citet{Fleming71,Fleming-Souga86}.
More recently, a fully probabilistic approach to the limit $v_n\to v$ was proposed by \citet{BaLaTa}.

\vspace{.2cm}
In the last decade, mostly motivated by the theory of mean field games and the control of McKean-Vlasov dynamics \cite{Huang2007,PLLcollege,MR2295621,MR3752669,MR3753660} and also by aspects of fluid mechanics and action minimizing paths \cite{Gang-Ngu09,Gangbo08}, there has been an intensive research activity around the Hamilton-Jacobi equation
\begin{equation}
\label{eq:HJ.McV}
	\begin{cases}
		\partial_t\cV(t,\mu) + \widetilde f(t, \partial_\mu \cV(t, \mu), \mu) = 0 \quad \text{in } [0,1) \times \cP_2(\R^\xdim)\\
		\cV(1, \mu) = \widetilde F(\mu) \quad \text{for } \mu \in \cP_2(\R^\xdim)
	\end{cases}
\end{equation}
where $\cP_p(\R^\xdim)$ is the space of probability measures on $\R^\xdim$ with finite $p$-th moment, which we equip with the Wasserstein metric $\cW_p$, and $\partial_\mu\cV$ denotes the Wasserstein gradient of $\cV$.

While this equation is less well understood than \eqref{eq:HJ}, impressive progress have been made recently.
The case 
\begin{equation}
\label{eq:examp.Hynd}
	\widetilde f(t, \partial_\mu\cV(t,\mu), \mu) := \frac1p\|\partial_\mu\cV(t,\mu)\|^p_{\L^p(\mu)} + \cU(\mu)
\end{equation}
for some energy potential $\cU:\cP_2(\R^\xdim) \to \R$ and $p>1$ has been particularly studied, starting with \citet{Gangbo08} who showed that when $\xdim =1, p=2$ and $\cU$ satisfies 
\begin{equation*}
	\cU(\mu) \le \alpha \cW_{\bar p}^{\bar p}(\mu, \nu) + \beta \quad \text{for some } \alpha,\beta \in \R \text{ and } \nu \in \cP_p(\R^\xdim)
\end{equation*}
with $\bar p$ the H\"older conjugate of $p$, and $\widetilde F$ and $\cU$ are continuous, then the following value function is a viscosity solution of Equation \eqref{eq:HJ.McV}:
\begin{equation}
\label{eq:value.Wass.space}
	\cV(s,\mu) = \sup_{\phi }\bigg(\widetilde F(\mu^\phi(1)) - \int_t^s\widetilde g\Big(u,\phi(u,\cdot), \mu^\phi(u) \Big)\diff u \bigg)
\end{equation}
where the supremum is over Borel--measurable maps $\phi:[0,1]\times \R^\xdim\to \R^\xdim$ and $\mu^\phi$ satisfies the continuity equation
\begin{equation*}
	\partial_t\mu^\phi + \nabla\cdot (\mu^\phi \phi) = 0, \quad \mu^\phi(s, \cdot) = \mu
\end{equation*}
in the sense of distributions and
\begin{equation*}
	\widetilde g(t, \partial_\mu\cV(t,\mu), \mu) := \frac1{\bar p}\|\partial_\mu\cV(t,\mu)\|^{\bar p}_{\L^{\bar p}(\mu)} - \cU(\mu).
\end{equation*}
This result was further extended by \citet{HyndJFA}.
See also \cite{Ambro-Feng,Gang-Tud19,Gang-Ngu09}.

Just as in the finite dimensional case \eqref{eq:HJ}, classical solutions of \eqref{eq:HJ.McV} are hard to expect in general (see however \citet{Gang-Swe15}), but conditions have been recently given (see \citet{ChassagneuxCrisanDelarue_Master} and \citet{carda15}) under which the following second order equation admits a unique classical solution:
\begin{align}\
	\label{eq:Wassers.PDE.n}
		\begin{cases}
			\partial_t\cV_n + \partial_x\cV_n\cdot b + \frac{1}{2n}\mathrm{Tr}\big[\partial_{xx}\cV_n\sigma\sigma^\top\big] + f\big(t, \partial_{x}\cV_n\sigma,x,\mu \big) \\
			\quad +\int_{\mathbb{R}^\xdim}\partial_\mu \cV_n(t, x, \mu)(a)b(t, a, \mu)\diff \mu(a) + \frac{1}{2n}\int_{\mathbb{R}^\xdim}\mathrm{Tr}\big[ \partial_a\partial_\mu \cV_n(t, x, \mu)(a)\sigma\sigma^\top(t,a,\mu)\big]\diff \mu(a)
			=0\\
			\cV_n(1, x,\mu) = F(x,\mu)\quad \text{for } (x,\mu) \in \R^\xdim \times \cP_2(\R^\xdim),
	\end{cases}
\end{align}
where the functions $b$ and $\cV_n$ are evaluated at $(t, x, \mu) \in [0,1]\times \R^\xdim \times \cP_2(\R^\xdim)$. (The link between the functions $F, f$ and $\widetilde F, \widetilde f$ is made clear below).
The main goal of this paper is to extend the vanishing viscosity results described above to PDEs on the space of probability measures.
Indeed, consider the following conditions:
\begin{itemize}
	\item[$(A1)$] The function $f:[0,1]\times \R^d\times \R^\xdim \times \cP_2(\R^\xdim)\to \R$ is such that
	$$f(t, z,x,\mu) = f_1(t, z) + f_2(t, x, \mu)$$ where the function $f_2(t,\cdot,\cdot)$ is Lipschitz--continuous and bounded, uniformly in $t\in [0,1]$, and $(x_1,\dots,x_N)\mapsto f_2(t, x_i, \sum_{j=1}^N\delta_{x_j})$ is continuous for all $N\ge1$ and all $i=1,\dots, N$.
	The function $f_1(t,\cdot)$ is convex, positive, satisfies $f_1(t,0)=0$ as well as the coercivity property 
	$$\lim_{\|z\|\to \infty}\frac{f_1(t, z)}{\|z\|}  = \infty$$ for all $t \in [0,1]$ and the integrability $\sup_{|z|\le r}f_1(t, z) \in \L^1([0,1],\diff t)$ for all $r \ge 0$.
	\item[$(A2)$] The function $b:[0,1]\times \R^\xdim\times \cP_2(\R^\xdim)\to \R^\xdim$ is $\ell_b$--Lipschitz--continuous, with $b$ of linear growth.
	That is,  
	\begin{equation*}
		\|b(t,x,\mu) - b(t,x^\prime,\mu^\prime)\| \le \ell_b\Big(\|x-x^\prime\| + \cW_2(\mu, \mu^\prime)\Big)  \text{ and } \|b(t,x,\mu)\|\le \ell_b\Big(1 + \|x\| + \big( \int_{\R^\xdim}\|y\|^2\mu(\diff y) \big)^{1/2} \Big)
	\end{equation*}
	for all $t \in [0,1]$, $x,x^\prime \in \mathbb{R}^\xdim$.
	\item[$(A3)$] $\sigma \in \R^{\xdim \times \bmdim}$ satisfies
	\begin{equation}
	\label{eq:sigma.elliptic}
	\langle y,\sigma\sigma^\top y\rangle>C_2 |y|^2\quad \text{for all }  y \in \mathbb{R}^\xdim \text{ for some } C_2>0.
	\end{equation}
\end{itemize}
The first main result of this paper links the solution $\cV_n$ of the second order equation \eqref{eq:Wassers.PDE.n} to (a general version of) the value function \eqref{eq:value.Wass.space}.
Therein, we denote by $g$ the convex conjugate of $f$
\begin{equation}
\label{eq:conv.conj.f}
 	g(t, q,x, \mu) := \sup_{z \in \R^\bmdim}\Big(q\cdot z - f(t, z, x, \mu)\Big).
 \end{equation} 
 Furthermore, we work on a filtered probability space $(\Omega, {\cal F},\P)$ carrying a standard $d$-dimensional Brownian motion $W$ and equipped with $({\cal F}_t)_{0\leq t\leq 1}$, the $\P$--completion of the filtration of $W$.
We use the notation
\begin{equation*}
	\widetilde F_s(\mu) := \int_{\R^\xdim}F(x,\mu)\mu(\diff x| \cF_s) \quad \text{and}\quad \widetilde g_s(t, \phi(t, \cdot), \mu) = \int_{\R^\xdim}g(t, \phi(t,x), x, \mu)\mu(\diff x| \cF_s)
\end{equation*}
for all $(t, \mu) \in [0,1]\times \cP_2(\R^\xdim)$ and $\mu(\cdot |\cF_s)$ the $\cF_s$--conditional distribution of $\mu$, with the tacit assumption that $F(\cdot,\mu), g(t,\phi(t,\cdot),\cdot, \mu) \in \L^2(\R^\xdim,\mu)$.
\begin{theorem}
\label{thm:PDE.convergence}
	Assume that the conditions $(A1)$, $(A2)$ and $(A3)$ are satisfied, and that $F$ is a bounded continuous function mapping $\R^\xdim\times \cP_2(\R^\xdim)$ to $\R$ and such that $(x_1,\dots, x_n)\mapsto F(x,\frac1n\sum_{i=1}^n\delta_{x_i})$ is continuous, for all $(x,n)\in \R^\xdim\times \mathbb{N}$.
	If a function $\cV_n:[0,1]\times \R^\xdim\times \cP_2(\R^\xdim) \to \R$ solves the PDE \eqref{eq:Wassers.PDE.n}, then for every $s \in [0,1]$ and $\xi\in \L^2(\R^\xdim,\cF_s)$ with absolutely continuous law $\mu$, it holds 
	\begin{equation*}
		\cV_n(s, \xi,\mu) \to \cV(s, \mu) := \esssup_{\phi}\bigg(\widetilde F_s(\mu^\phi(1)) - \int_s^1\widetilde g_s\Big(u, \phi(u, \cdot,\mu^\phi(u)), \mu^\phi(u)\Big)\diff u \bigg)\quad \P\text{-a.s.}
	\end{equation*}
	where the supremum is over Borel--measurable maps $\phi:[0,1]\times \R^\xdim\times \cP_2(\R^\xdim)\to \R^\xdim$ such that $\phi(t,\cdot, \mu)\in \L^2(\R^\bmdim, \mu^\phi)$ with $\mu^\phi$ satisfying the continuity equation
	\begin{equation*}
		\partial_t\mu^\phi + \nabla\cdot \big( \mu^\phi b(t, \cdot, \mu) + \mu^\phi \sigma\phi(t, \cdot, \mu^\phi) \big) = 0, \quad \mu^\phi(s, \cdot) = \mu
	\end{equation*}
	in the sense of distributions.

	If the function $b$ does not depend on $\mu$, i.e. $b(t, x, \mu) = b(t,x)$, then it holds
	\begin{equation*}
		\cV_n(s, \xi,\mu) \to \cV(s, \mu) := \esssup_{\phi}\bigg(\widetilde F_s(\mu^\phi(1)) - \int_s^1\widetilde g_s\Big(u, \phi(u, \cdot), \mu^\phi(u)\Big)\diff u \bigg)\quad \P\text{-a.s.}
	\end{equation*}
	where the supremum is over Borel--measurable maps $\phi:[0,1]\times \R^\xdim\to \R^\xdim$ such that $\phi(t,\cdot) \in \L^2(\R^\bmdim,\mu^\phi)$ and $\mu^\phi$ satisfies 
	\begin{equation*}
		\partial_t\mu^\phi + \nabla\cdot \big( \mu^\phi b(t, \cdot) + \mu^\phi \sigma\phi(t, \cdot) \big) = 0, \quad \mu^\phi(s, \cdot) = \mu.
	\end{equation*}
\end{theorem}
The proof is given in the subsection \ref{sec:PDE.proof}.
The case \eqref{eq:examp.Hynd}-\eqref{eq:value.Wass.space} is obtained with the specifications
\begin{equation*}
	b=0\quad \text{and}\quad f(t, z, x, \mu) = \frac1p\|z\|^p + \cU(\mu).
\end{equation*}
Observe that the limiting function $\cV(t, \mu)$ of $\cV_n(t, \xi, \mu)$ depends on $\xi$ through its law, since $\mu = \text{law}(\xi)$.

The above convergence result seems to be the first of its kind for equations on the Wasserstein space.
Note however the work of \citet{Cecc-Del20}.
In this work, the authors (among other things) show the convergence of the master equation of a finite state mean field game with common noise converges to (the gradient of) the mean field control problem without common noise, showing for instance how vanishing common noise allows to select a particular solution of the mean field game.
See also \cite{Del-Tcheu20} for a result along the same lines.
In the present case, the limiting control problem is one with deterministic states (first order Fokker-Planck equations).

\vspace{.2cm}
The proof of Theorem \ref{thm:PDE.convergence} is based on fully probabilistic arguments in the spirit of those developed by \citet{BaLaTa} for the case of equations on Euclidean spaces.
In fact, the cornerstone of our approach is a variational representation formula for non-exponential functions of McKean-Vlasov diffusions which generalizes the celebrated Gibbs variational principle first discovered by \citet{Fl78,Boue-Dup} in the case of Wiener process.
This variational principle is a stochastic control representation of the cumulant generating function of Wiener process.
See also \cite{Bo00,Lehec}.
En route to proving Theorem \ref{thm:PDE.convergence}, we derive a novel large deviation theorem in the form of a non-exponential Freidlin-Wentzell theorem in its Laplace principle form, see Theorem \ref{thm:abstract.schilder}.
A special case of this theorem will lead to the Freidlin-Wentzell theorem for McKean-Vlasov diffusions.
In fact, we consider the following condition:
\begin{itemize}
\item[$(A3)'$] The function $\sigma:[0,1]\times \R^\xdim\times \cP_2(\R^\xdim)\to \R^{\xdim\times \bmdim}$ is $\ell_b$--Lipschitz--continuous, and bounded.
	That is,  
	\begin{equation*}
		\|\sigma(t,x,\mu) - \sigma(t,x^\prime,\mu^\prime)\| \le \ell_b\Big(\|x-x^\prime\| + \cW_2(\mu, \mu^\prime)\Big)  \text{ and } \|\sigma(t,x,\mu)\|\le \ell_b
	\end{equation*}
	for all $t \in [0,1]$, $x,x^\prime \in \mathbb{R}^\xdim$, $\ell_b>0$, and 
	$$\langle y,\sigma(t,x,\mu)\sigma^\top(t',x,\mu)y\rangle>C_2 |y|^2\quad \text{for all }  (t,,\mu,x,y) \in [0,1]\times\cP_2(\R^\xdim)\times (\mathbb{R}^\xdim)^2 \text{ for some } C_2>0.
	$$
\end{itemize}
In the statement of the result, we will use the space $\mathcal{H}$ defined as
\begin{equation}
\label{eq:def.calH}
 	\mathcal{H}:=\Big\{\varphi:[0,1]\to \R^\bmdim: \text{Borel--measurable and } \int_0^1\|\varphi(t)\|^2\diff t<\infty \Big\}.
 \end{equation}\
\begin{corollary}[Freidlin-Wentzell Theorem]
\label{cor:F-WThm}
	Assume that the conditions $(A2)$ and $(A3)'$ are satisfied.
	Given $x \in \R^\xdim$ and $n\ge1$, let $X_n$ solve the SDE
	\begin{equation}
	\label{eq:McKV.SDE.n}
		\begin{cases}
		dX_n(t) = b\big(t, X_n(t), \mu_n(t)\big)\diff t + \frac{1}{\sqrt{n}}\sigma\big(t, X_n(t), \mu_n(t)\big)\diff W(t),\\
		 \quad X_n(0) = x,\,\, \mu_n(t) = \text{law}(X_n(t)).
		\end{cases}
	\end{equation}
	For every bounded continuous functions $F :\mathbb{R}^m\times \cP_2(\mathbb{R}^m) \to \R$, it holds
	\begin{equation*}
		\lim_{n\to \infty}\frac1n\log\E[e^{nF(X_n(1), \mu_n(1))}] = \sup_{ \varphi \in \mathcal{H}}\Big(F(\Phi^\varphi(1),\delta_{\Phi^\varphi(1)}) - I(\Phi) \Big)
	\end{equation*}
	with\footnote{As usual, we adopt the convention $\inf{\emptyset}: = +\infty$.}
		\begin{equation*}
		I(\Phi) := \inf\Big\{\frac12\int_0^1\|\varphi(t)\|^2\diff t: \varphi \in \mathcal{H}\,\text{s.t. } \Phi(t) = x + \int_0^tb(u, \Phi(u), \delta_{\Phi(u)})+ \sigma(u, \Phi(u), \delta_{\Phi(u)})\varphi(u)\diff u  \Big\}.
	\end{equation*}
\end{corollary}
As stated above, this result is an extension of the celebrated Freidlin-Wentzell Theorem to the case of McKean-Vlasov diffusion, and is thus interesting in its own right.
The work of \citet{Herr-Imkeller-Pei08} probably gives one of the first results in this direction.
In fact, these authors prove a small noise large deviation result for a McKean-Vlasov diffusion with constant diffusion term $\sigma$ and a specific drift term with linear dependence in the measure argument.
This work was further extended (with simpler proofs) by \citet{Tug16}.
The first general Freidlin-Wentzell Theorem for McKean-Vlasov diffusions is due to \citet{dReis-Sal-Tug18} who give a similar result to ours.
Their results hold for functionals $F$ on the path space, but the assumptions made on the functions $b$ and $\sigma$ are more restrictive than ours.
In addition to the conditions imposed on the coefficients $b$ and $\sigma$, Corollary \ref{cor:F-WThm} (or actually its more general version Theorem \ref{thm:abstract.schilder}) is also interesting due to its method of proof.
Indeed, the argument for this result builds on the ``weak compactness'' approach to large deviations developed by \citet{dupuis-ellis} and applied to the standard Freidlin-Wentzell theorem by \citet{Boue-Dup} and \citet{BaLaTa}.

In this work we further consider a consequence of our main representation result to showcase its relevance beyond the vanishing viscosity problem.
We will take advantage of our extended Gibbs variational principle to the derivation of a functional inequality for linear McKean-Vlasov equations.
Following the original idea of \citet{Bo00}, we prove that solutions of McKean-Vlasov equations satisfy the Pr\'ekopa-Leindler inequality (a reverse form of H\"older's inequality). 
This functional inequality first proved by \citet{Prekopa71} and \citet{Leindler72} to study problems in linear programing has turned out over the years to have fundamental applications in analysis, geometry and probability theory, We refer the interested reader to \cite{Gardner02} for an extensive overview and further applications.
The modest result provided here is a simple observation that should be further extended in the future.

Most of the rest of the paper is dedicated to the proof of Theorem \ref{thm:PDE.convergence}.
We will start by introducing a convex functional generalizing the log-moment generating function and study its variational representation when applied to (functions of) McKean-Vlasov diffusions.
This result will allow to derive a new proof of the Freidlin-Wentzell theorem stated in Corollary \ref{cor:F-WThm}.
An extension of this result along with a version of the Feynman-Kac formula will allow to easily  conclude the proof of Theorem \ref{thm:PDE.convergence}.
In the final section, we discuss an application of to functional inequalities.

\section{A Variational representation}

In this section we prove a new variational representation for functional of McKean-Vlasov dynamics.
This representation will play a crucial role in the ensuing proofs of our main results.
We start by presenting the probabilistic setting and some notation.
\subsection{Preliminaries}
\label{sec setting}
Recall that we work on a filtered probability space $(\Omega, {\cal F},\P)$ carrying a standard $d$-dimensional Brownian motion $W$ and equipped with the $\P$--completion of the filtration of $W$ denoted $\F:=({\cal F}_t)_{0\leq t\leq 1}$.

Consider the space
\begin{equation*}
	\mathcal{L}:=\left\{q: \Omega\times [0,1] \to \mathbb{R}^d; q\text{ is progressive, and }\int_{0}^{T}\|q(t)\|^2\diff t <+\infty \,\, \P\text{-a.s.}\right\}.
\end{equation*}
As usual, we identify random variables that are equal $\P$-a.s. and processes that are indistinguishable.
We will denote by $\cL^\infty$ the elements of $\cL$ which are bounded, and by $\cL_k^\infty$ the elements of $\cL^\infty$ bounded by the positive number $k$.
Furthermore, given a non-empty metric space $E$, a sigma algebra $\Gc$ and a filtration $\G$, we will use the following norms:

\medskip
$\bullet$ For any $p\in[1,\infty]$, $\L^p(E,\Gc)$ is the space of $E$-valued, $\Gc$-measurable random variables $R$ such that 
\[
	\|R\|_{\L^p(E,\Gc)}:=\Big(\E\big[\|R\|_E^p\big]\Big)^{\frac1p}<\infty,\; \text{when}\; p<\infty,\; \|R\|_{\L^\infty(E,\Gc)}:=\inf\big\{\ell\geq0:\|R\|_E\leq \ell,\; \P\text{\rm --a.s.}\big\}<\infty.
\]

\medskip
$\bullet$ For any $p\in[1,\infty)$, $\H^p(E,\G)$ is the space of $E$-valued, $\G$-predictable processes $Z$ such that 
\begin{equation*}
	\|Z\|_{\H^p(E,\G)}^p:=\EE\bigg[\bigg(\int_0^T\|Z_s\|_E^2\mathrm{d}s\bigg)^{p/2}\bigg]<\infty.
\end{equation*}

$\bullet$ For any $p\in[1,\infty]$, $\S^p(E,\G)$ is the space of $E$-valued, continuous, $\G$-adapted processes $Y$ such that 
\begin{equation*}
	\|Y\|_{\S^p(E,\G)}:=\bigg(\EE\bigg[\sup_{t\in[0,T]}\|Y_t\|_E^p\bigg]\bigg)^{\frac1p}<\infty,\;  \text{when}\; p<\infty,\; \|Y\|_{\S^\infty(E,\G)}:=\bigg\|\sup_{t\in[0,T]}\|Y_t\|_E\bigg\|_{\L^\infty(E,\Gc_T)}<\infty.
\end{equation*}
Furthermore, we will always denote by $\|f\|_\infty$ the smallest upper bound of a function $f$, irrespective of the space on which it is defined.
We will also denote by $\cC_b(E)$, the space of continuous functions on $E$, and by $\cC_b(E)$ the bounded elements of $\cC(E)$.

Let ${\cal Q}$ be the set of probability measures absolutely continuous w.r.t. $\P$.
It is well-known that for every $\Q \in {\cal Q}$, there is a unique process  $q^\Q \in {\cal L}$ such that
\begin{equation*}
	\frac{\diff \Q}{\diff \P} = \exp\bigg(\int_0^1q^\Q(t)\diff W(t)-\int_0^1\frac{1}{2}\|q^\Q(t)\|^2\diff t \bigg).
\end{equation*}
Given a function $f:[0,1]\times \R^\bmdim\times \R^\xdim\times \cP_2(\R^\xdim)\to \R$ that is convex it its second argument and we denote by $g$ its convex conjugate defined in \eqref{eq:conv.conj.f}.
For every $s \in [0,1]$ we define the functional $\alpha^g_s: {\cal Q}\to \L^0\big(\mathbb{R}\cup\{+\infty\} , {\cal F}_s\big)$ and its conjugate $\rho^g_s:\L^0(\R,\cF_T)\to \L^0\big(\mathbb{R}\cup\{+\infty\} , {\cal F}_s\big)$, respectively given by
\begin{equation*}
	\alpha^g_s(\Q):= \E^\Q\left[\int_s^1g\big(t,q^\Q(t),X(t),\mu(t) \big)\diff t\mid{\cal F}_s\right] 
\end{equation*}
and
\begin{equation*}
	 \rho^g_s(G):=\esssup_{\Q \in  {\cal Q}} \left(\E^\Q[G\mid{\cal F}_s]-\alpha^g_s(\Q)\right),
\end{equation*}
where $X$ is the (strong) solution of the McKean-Vlasov equation
\begin{equation}
	\diff X(t) =   b(t, X(t), \mu(t) )\diff t + \sigma(t, X(t), \mu(t) )\diff W(t),\quad X(0) = \xi, \,\, \mu(t) =\text{law}(X(t)).
\end{equation} 
It is interesting to notice that, when $f(t,z,x,\mu) :=  \frac12\|z\|^2$, then the conjugate $g$ is again the quadratic function $g(t,q,x,\mu) = \frac12\|q\|^2$.
In this case, the functional $\alpha^g_s$ is the (conditional) Kullback-Leibler divergence and, by Gibbs' variational principle, $\rho^g_s$ is nothing other than the (conditional) cumulant moment generating functional, i.e.
\begin{equation*}
 	\alpha^g_s(\Q) = \E^{\Q}\Big[\frac{\diff \Q}{\diff \P}\log\Big(\frac{\diff \Q}{\diff \P}\Big) \Big| \cF_s\Big]\quad \text{and}\quad \rho^g_s(G) = \log\big( \E[e^G| \cF_s] \big).
\end{equation*} 
This is the quintessential example in our analysis that the reader should always keep in mind.

The main result of this section is a variational representation of the functional $\rho^g_s$.
We will prove the following:
\begin{theorem}
\label{thm:var.rep.McKV}
	Assume that the conditions $(A1)$, $(A2)$ and $(A3)'$ are satisfied.
	Let $s \in [0,1]$ and $\xi \in \L^2(\R^\xdim,\cF_s)$.
	For every lower semicontinuous function $F:\R^\xdim\times \cP_2(\R^\xdim)\to [-c,\infty]$ for some $c\ge 0$ such that $(x_1,\dots, x_n)\mapsto F(x,\frac1n\sum_{i=1}^n\delta_{x_i})$ is continuous for all $n\in \mathbb{N}$, it holds 
	\begin{align}
	\notag
		& \rho^{g}_s\big( F(X^{s,\xi}(1), \mu^{s,\xi}(1)) \big)\\\label{eq:var.rep.abstract.McKV}
		&\qquad\qquad  = \esssup_{q\in \cL}\E\bigg[F\big( X^{s,\xi,q}(1), \mu^{s,\xi,q}(1) \big) - \int_s^1 g\Big(t,q(t),X^{s,\xi,q}(t), \mu^{s,\xi,q}(t) \Big)\diff t\bigg|\cF_s \bigg],
	\end{align}
	where $X^{s,\xi,q}$ is the strong solution of the (controlled) McKean-Vlasov equation
	\begin{align}
	\notag
		X^{s,\xi,q}(t) &= \xi + \int_s^tb\big(u, X^{s,\xi,q}(u), \mu^{s,\xi,q}(u)\big) + \sigma\big(u, X^{s,\xi,q}(u), \mu^{s,\xi,q}(u)\big)q(u)\diff u\\
		\label{eq:contr.MckV}
		&\quad + \int_s^t\sigma\big(u, X^{s,\xi,q}(u), \mu^{s,\xi,q}(u)\big)\diff W(u),\quad \mu^{s,\xi,q}(u) = \mathrm{law}(X^{s,\xi,q}(u))
	\end{align}
	with the convention $X^{s,\xi}:= X^{s,\xi,0}$.
\end{theorem}
\begin{remark}
	Let us make the following observations:
	\begin{itemize}
		\item The continuity condition on the map $(x_1,\dots, x_n)\mapsto F(x,\frac1n\sum_{i=1}^n\delta_{x_i})$ is clearly guaranteed, for instance when $F(x,\cdot)$ is Lipschitz--continuous with respect to the Wasserstein distance of any order.
		\item Well-posedness of the McKean-Vlasov equation \eqref{eq:contr.MckV} is standard under our assumptions.
		See e.g.\ \cite{MR3752669}.
\end{itemize}
\end{remark}

\vspace{.2cm}

The representation Theorem \ref{thm:var.rep.McKV} was first discovered for the case $g(t,q,x,\mu):=\frac12\|q\|^2$, and appeared in \cite{Bo00,Boue-Dup,Fl78}.
These results were either proved for Brownian motion (i.e.\ when $b=0$ and $\sigma$ is a constant) or for classical SDEs (i.e. when $b$ and $\sigma$ do not depend on the law of the unknown).
For functions $F$ of Brownian motion and $g$ depends only on $(t,q)$, Theorem \ref{thm:var.rep.McKV} was derived in \cite{BaLaTa}.
All these variants of the variational representation \eqref{eq:var.rep.abstract.McKV} have had crucial consequences in surprisingly different areas.
For instance, in large deviations theory \cite{Boue-Dup,BaLaTa} in convex geometry and functional inequalities \cite{Lehec,Bo00,vanHandel17} or in the study of the Schr\"odinger problem in optimal transportation theory \cite{BaLaTa}.
The extension of the current paper concerns McKean-Vlasov SDEs and further allows $F$ and $g$ to depend on the law of the SDE in question.
These extensions force the introduction of a proof of the representation that is very different from those proposed in all the aforementioned references.
Our proof is tailor-made for the case where $F$ depends on the terminal value of $X$, (not its entire path as in \cite{Boue-Dup,BaLaTa}).
But this is sufficient for the main motivation of the present article. 
In a first step, we will derive the representation for the case of classical SDEs.
This result seems to be interesting on its own right as it extends the representation of \cite{BaLaTa} to the case of functions of SDEs. 
Next, we will use standard propagation of chaos arguments to finish the proof of Theorem \ref{thm:var.rep.McKV}.

\subsection{Proof of the variational representation}
In preparation of the proof of Theorem \ref{thm:var.rep.McKV}, let $\ell\in \N$ be given and consider two functions $B:[0,1]\times \R^\ell\to \R^\ell$ and $\Sigma :[0,1]\times \R^\ell\to \R^{\ell\times d}$.
For every $q \in \cL$, denote by $X^{s,\xi,q}$ the solution of the SDE
\begin{align}
	\label{eq:controlled.SDE}
		X^{s,\xi,q}(t) = \xi  + \int_s^t\Big\{ B\big( u,X^{s,\xi,q}(u) \big) + \Sigma\big( u,X^{s,\xi,q}(u) \big)q(u)\Big\}\diff u
		\ + \int_s^t\Sigma\big( u,X^{s,\xi,q}(u) \big)\diff W(u)
\end{align}
for some  $\xi \in \L^2(\R^\ell,\cF_s)$, again with the convention $X^{t,\xi} := X^{t, \xi, 0}$.
Consider the following conditions:
\begin{itemize}
	\item[{(A1)'}] The functions $B:[0,1]\times \R^\ell\to \R^\ell$ and $\Sigma:[0,1]\times \R^\ell\to \R^{\ell\times d}$ are $\ell_B$--Lipschitz--continuous and bounded (where $\ell_B$ does not depend on $\ell$).
	That is,  
	\begin{equation*}
		|B(t,x) - B(t,x^\prime)| + |\Sigma(t,x) - \Sigma(t,x^\prime)|\le \ell_B|x-x^\prime| \quad \text{and } \|\Sigma\|_\infty +\|B\|_\infty\le \ell_B
	\end{equation*}
	for all $t \in [0,1]$, $x,x^\prime \in \mathbb{R}^\ell$,
	and in addition 
	$$\langle y,\Sigma(t,x)\Sigma^\top(t',x)y\rangle>C_2 |y|^2\quad \text{for all }  x,y \in \mathbb{R}^\ell \text{ for some } C_2>0.
	$$
	\item[{(A2)'}] The function $g:[0,1]\times \R^d\times \R^\ell \to \R$ is such that $g(t,q,\cdot)$ is continuous for each $(t,q)\in [0,1]\times \R^\bmdim$, and the function $g(t,\cdot,x)$ is convex, lower semicontinuous, positive satisfies $g(t,0,x)=0$ as well as the coercivity property $\lim_{\|q\|\to \infty}\inf_{(t,x)\in[0,1]\times\R^\xdim}\frac{g_1(t, q,x)}{\|q\|}  = \infty$ and the integrability $\sup_{\|(q,x)\|\le r}g(t, q,x) \in \L^1([0,1],\diff t)$ for all $r \ge 0$.
\end{itemize}
\begin{proposition}
\label{prop:var.rep}
	If the conditions $(A1)'$ and $(A2)'$ are satisfied, then for every bounded, lower semicontinuous function $F:\R^\ell\to \R$, 
	\begin{align}
	\label{eq:var.rep.abstract.}
		& \rho^{g}_s\big( F(X^{s,\xi}(1)) \big)  = \esssup_{q\in \cL^\infty}\E\left[F\big( X^{s,\xi,q}(1) \big) - \int_s^1 g\Big(t,q(t),X^{s,\xi,q}(t) \Big)\diff t\bigg|\cF_s \right].
	\end{align}
\end{proposition}
An essential element in the proof of this proposition is the link between the operator $\rho^g$, solutions of backward stochastic differential equation (BSDE) and semi-linear parabolic PDEs.
We refer the unfamiliar reader for instance to the articles \cite{Pardoux-Peng92,karoui01,DHK1101,tarpodual,kobylanski01,Delbaen11} or to the recent monograph \cite{zhangbook}.
For the reader's convenience, we summarize the results that will be needed here in the following lemma:
\begin{lemma}
\label{lem:BSDE.misc}
	Let $F\in \L^2(\R, \cF_T)$ be a given random variable and $f:[0,1]\times\Omega\times \R^d\to \R$ be a given function.
	\begin{itemize}
		\item[(i)]	
		If $f$ is Lipschitz--continuous in its last argument, $f(\cdot, y,z)$ is progressively measurable and satisfies $\E\Big[\int_0^T|f(t,0)|^2\diff t\Big]<\infty$, then there is a unique pair $(Y, Z)\in \mathbb{H}^2(\mathbb{R},\mathbb{F})\times \mathbb{H}^2(\R^d,\FF)$ solving the BSDE with terminal condition $F$ and generator $f$, i.e.
	\begin{equation}
	\label{eq:lem.bsde}
		Y(t) = F + \int_t^1f(u, Z(u))\diff u - \int_t^1Z(u)\diff W(u).
	\end{equation}
		\item[(ii)]
		If the random variables $F$ and $f(t,z)$ are of the form $F\equiv F(X^{s,x}(1))$ and $f(t,z)\equiv f(t, z, X^{s,x}(t))$, then $ v(s, x):=Y(s)$ is a viscosity solution of the PDE
		\begin{equation*}
			\begin{cases}
				\partial_tv + \partial_xv\cdot B + \frac{1}{2}\mathrm{Tr}\big[\partial_{xx}v\Sigma\Sigma^\top\big] + f(t, \partial_{x}v\Sigma,x) = 0 \quad \text{in } [0,1)\times \R^\ell,\\
				v(1, x) = F(x)\quad \text{for } x\in \R^\ell.
			\end{cases}	
		\end{equation*}
		\item[(iii)] Assume that \eqref{eq:lem.bsde} admits a solution $(Y,Z)$ and that $F$ is bounded.
		If $f$ satisfies $|f(t,z)|\le C(1 + \|z\|^2)$ for some $C>0$, is bounded from below, is convex and lower semicontinuous and in $z$, then it holds
		\begin{equation*}
			Y_s = \rho^g_s(F) \quad \PP\text{-a.s.}
		\end{equation*}
		where $g$ is the convex conjugate of $f$.
	\end{itemize}
\end{lemma}
\begin{proof}
	The statement $(i)$ is \cite[Theorem 2.1]{karoui01}.
	The statement $(ii)$ is \cite[Theorem 4.2]{karoui01} and $(iii)$ is \cite[Proposition 3.1 $\&$ Remark 3.6]{tarpodual}.
\end{proof}

\begin{proof}[Proof of Proposition \ref{prop:var.rep}]
	In preparation for the proof, let us introduce the concept of minimal supersolutions of backward stochastic differential equations.
	Henceforth, denote by $f$ the function defined on $[0,1]\times \mathbb{R}^d\times \R^\ell $ by
	\begin{equation*}
		f(t,z,x) := \sup_{q\in \mathbb{R}^d}\Big( q\cdot z - g(t, q, x) \Big).
	\end{equation*}
	Following \cite{DHK1101} we call a pair $(Y,Z)$ with $Y$ a real--valued c\'adl\'ag and adapted process and $Z\in {\cal L}$ a supersolution of the backward stochastic differential equation 
	\begin{equation}\textstyle
	\label{eq:bsde}
		\diff Y(t) = -f\big(t,Z(t), X^{s,\xi}(t) \big)\diff t + Z(t)\cdot\diff W(t), \quad Y(1) = F\big( X^{s,\xi}(1) \big)
	\end{equation}
	with terminal condition $F(X^{s,\xi}(1))$ and generator $f$ if it satisfies
	\begin{equation}\textstyle
		\begin{cases}
 			Y(s)-\int_{s}^{t}f(u,Z(u), X^{s,\xi}(u))\diff u+\int_{s}^{t}Z(u)\cdot \diff W(u)\geq Y(t), \quad \text{for every} \quad 0\leq s\leq t\leq 1\\
			\displaystyle Y(1)\geq F(X^{s,\xi}(1))
		\end{cases}
		\label{eq:supersolutions}
	\end{equation}
	and $\int Z\cdot\diff W$ is a supermartingale.
	A supersolution $(\bar{Y},\bar{Z})$ of \eqref{eq:bsde} is said to be minimal if $\bar{Y}(t)\leq Y(t)$ for all $t \in [0,1]$ and for every other supersolution $(Y,Z)$.
	Introducing minimal supersolutions is needed here because very weak growth assumptions are made on $f$.
	It is well-known, (see e.g. \cite{Delbaen11}) that unless $f$ is at most of quadratic growth, the BSDE is not well-posed.
	But a minimal supersolution as defined above can still be defined. 
	In fact,
	since for almost every $(t,x)\in [0,1]\times\R^\ell$, $f(t,\cdot,x)$ is lower semicontinuous, convex, and positive, by \cite[Theorem 4.5]{tarpodual}, there exists a unique $\bar Z \in {\cal L}$ such that $(\bar Y,\bar Z)$ is the minimal supersolution of Equation \eqref{eq:bsde} and more importantly for our purpose here, it holds 
	$$\bar Y(t) = \rho^g_t\big(F(X^{s,\xi}(1) \big)\quad \P\text{-a.s.}$$
	With this representation at hand, we can split the rest of the proof in the following two steps.

\vspace{.2cm}

	\emph{Step 1: The lower bound.}  

	\emph{Step 1a: The Lipschitz continuous case.} We first assume that there is $K\ge 0$ such that
	\begin{equation*}
		\|F(x) - F(x')\|\le K\|x-x'\|  \quad \text{and} \quad \|F\|\le K
	\end{equation*}
	for all  $(t, x,x^\prime) \in [0,1]\times (\mathbb{R}^\ell)^2$.
	For every $n \in \mathbb{N}$, define the truncated function
	\begin{equation}
	\label{eq:trunc.f}
		f_n(t,z,x):=\sup_{ \|q\|\le n}\Big(q \cdot z - g(t,q,x) \Big).
	\end{equation}
	Under $(A1)'$, and for  $q \in {\cal L}_n$, the equation \eqref{eq:controlled.SDE} has Lipschitz--continuous coefficients and therefore admits a unique strong solution $X^{s,\xi,q}$.
	Moreover, since $f_n$ is Lipschitz--continuous, the BSDE with terminal condition $F\big( X^{s,\xi}(1) \big)$ and generator $f_n$ admits a unique square integrable solution $(Y_n,Z_n)\in \H^2(\R, \F)\times \H^2(\R^d, \F)$, see Lemma \ref{lem:BSDE.misc}. 
	In addition, it follows again from Lemma \ref{lem:BSDE.misc} that $Y_n(t)= \rho^{g_n}_t\big( F(X^{s,\xi}(1)) \big)$ with $g_n$ the convex conjugate of $f_n$ given by
	\begin{equation}
	\label{eq:conjug.trunc}
		g_n\big(t,q,x\big):=\sup_{z\in \mathbb{R}^d}\Big( q\cdot z - f_n(t,z,x)\Big).
	\end{equation}
	On the other hand, it follows by \cite[Theorem 8.4]{Grandall-Koc-Sci} that there is a function $v_n$ solving PDE 
	\begin{align}\
	\label{eq:PDE.bsde.n}
		\begin{cases}
			\partial_tv_n + \partial_xv_n\cdot B + \frac{1}{2}\mathrm{Tr}\big[\partial_{xx}v_n\Sigma\Sigma^\top\big] + f_n(t, \partial_{x}v_n\Sigma,x) \quad \text{in } [0,1)\times \R^\ell	=0,\\
			v_n(1, x) = F(x)\quad \text{for } x\in \R^\ell
		\end{cases}
	\end{align}
	in the Sobolev sense.
	By Lemma \ref{lem:BSDE.misc} it holds that $Y_n(s) = v_n(s,\xi)$.
	Let $ \Q\in {\cal Q}$ be such that $q\equiv q^\Q$ is bounded.
	Applying It\^o's formula for Sobolev functions (see \cite[Theorem 10.1]{KrylovBook}) to the process $v^n\big(t, X^{s,\xi,q}(t) \big)$, one has
	\begin{align*}\
		\diff v_n\big(t, X^{s,\xi,q}(t)\big) 
		& = \Big(  \partial_tv_n+\partial_xv_n\cdot b + \partial_xv_n \Sigma\cdot q(t) + \frac{1}{2} \mathrm{Tr}\big[\partial_{xx}v_n\Sigma\Sigma^\top\big]   \Big)(t, X^{s,\xi,q}(t))\diff t\\
		&\qquad + \partial_xv_n\Sigma\big(t, X^{s,\xi,q}(t) \big)\cdot \diff W(t).
	\end{align*}
	Therefore, since $v_n$ solves \eqref{eq:PDE.bsde.n}, it holds
	\begin{align*}
		F\big(X^{s,\xi,q}(1)\big) - v_n\big(s,\xi\big)
		& = \int_s^1 \partial_xv_n \Sigma\cdot q(t)-f_n\big(t,\partial_xv_n\Sigma,X^{s,\xi,q}(t) \big) \diff t + \int_s^1 \partial_xv_n\Sigma\cdot\diff W(t),
	\end{align*}
	where in the above equation, when no argument is given, the functions are evaluated at $(t, X^{s,\xi,q}(t))$.
	Taking conditional expectation above, and letting $n$ be large enough that $\|q\|\le n$, one obtains by definition of $g_n$ that 
	\begin{align}
	\label{eq:upper.bound.lips.n}
		\rho^{g_n}_s\big( F(X^{s,\xi}(1)\big) & = v_n(s,\xi)
		 \ge \E\left[F\big( X^{s,\xi,q}(1) \big) - \int_s^1 g_n\Big(t,q(t),X^{s,\xi,q}(t) \Big)\diff t\bigg|\cF_s \right].
	\end{align}
	Since $f_n\le f$, it follows that $g\le g_n$.
	Thus, by \cite[Proposition 3.3]{DHK1101},
	\begin{equation*}
		\rho^{g_n}_s\big( F(X^{s,\xi}(1) \big)\le \rho^g_s\big( F(X^{s,\xi}(1) \big) \quad\P\text{--a.s.}
	 \end{equation*} 
	Therefore, taking the limit in \eqref{eq:upper.bound.lips.n} it follows by monotone convergence that
	\begin{align}
	\label{eq:upper.bound.lips}
		\rho^{g}_s\big( F(X^{s,\xi}(1)\big) 
		 \ge \E\bigg[F\big( X^{s,\xi,q}(1) \big) - \int_s^1 g\Big(u,q(u),X^{s,\xi,q}(u) \Big)\diff u \bigg|\cF_s \bigg]\quad \text{for all } q \in \cL^\infty.
	\end{align}

	\emph{Step 1b}:
	Now, if $F$ is lower semicontinuous and bounded, say by $k$, then denote by 
	$$F_m(x):=\inf_{y\in \mathbb{R}^d}\big( F(y)+ m|x - y| \big)$$
	the so-called Pasch-Hausdorff envelope of $F$.
	For each $m$, the function $F_m$ is bounded and Lipschitz continuous; and $F_m\uparrow F$ pointwise, see e.g. \cite[Example. 9.11]{rockafellar02} for details.
	It follows from \emph{Step 1a} that for every $m \in \mathbb{N}$ and $q \in {\cal L}$, the bound \eqref{eq:upper.bound.lips} holds.
	Moreover, for all $m\ge1$, $\rho^{g}_s\big( F_m(X^{s,\xi}(1) \big)\le \rho^g_s\big( F(X^{s,\xi}(1) \big)$ $\P$-a.s., see \cite[Proposition 3.3]{DHK1101}.
	Therefore, we deduce that \eqref{eq:upper.bound.lips} holds for $F$ and $g$.

	\vspace{.2cm}
	
	\emph{Step 2:}
	Assume that $F$ is bounded, say by $K\ge0$.
	As above, we start with the identity
	$$\rho^{g_n}_t\big( F(X^{s,\xi}(1) \big) = v_n\big(t, X^{s,\xi}(t) \big)$$
	 for every $n \in \mathbb{N}$, 
	 where $v_n$ is the solution of the PDE \eqref{eq:PDE.bsde.n}.
	For each $(t,x)\in [0,T]\times\R^\ell$, let $ q(t,x)\in \mathbb{R}^d$ be in the subgradient of $f_n(t,\cdot, x)$ at $\partial_xv_n(t,x)\Sigma$.
	That is,  
	\begin{align}
		&f_n\big(t, (\partial_xv_n\Sigma)(t, x),x \big) = (\partial_xv_n\Sigma\cdot q)(t,x) - g_n\big(t, q(t,x),x\big).
		\label{eq:subgrad}
	\end{align}
	Since $f_n$ is Lipschitz--continuous in $z$ (with Lipschitz constant $n$), the function $q$ is bounded, see e.g. \cite[Lemma 3.2]{tarpodual}, and can be chosen Borel--measurable.
	The latter claim following by application of  \cite[Theorem 14.56]{rockafellar02} and the measurable selection theorem \cite[Corollary 1C]{rockafellar03},  since $\partial_xv_n\Sigma$ is Borel measurable on $[0,1]\times \R^\ell$, and for each $(t,x)$ the subgradient of $f_n(t,\cdot,x)$ at $\partial_xv_n\Sigma(t,x)$ is non-empty.
	Thus, the following multidimensional SDE with bounded, Borel--measurable drift and Lipschitz--continuous volatility admits a strong solution:
	$$
		dX(t) = B\big( t,X(t) \big) + q(t,X(t))\Sigma\big(t,X(t) \big)\diff t + \Sigma\big(t,X(t)\big)\diff W(t)
	$$
	admits a strong solution, see \cite[Theorem 1]{Verete_USSR}.
	With the notation of \eqref{eq:controlled.SDE}, we have $X^{s,\xi,q(\cdot, X_\cdot)} = X$ and recall that by construction that $q(\cdot,X_\cdot)\in \cL^\infty$.
	Applying It\^o's formula this time to $v_n(t,X(t))$, one gets
	\begin{align*}
		v_n\big( 1,X(1)\big) - v_n\big(s,\xi \big)
		&=  \int_s^1\left(  \partial_tv_n + \partial_xv_n\cdot B + \frac{1}{2}\mathrm{Tr}[\partial_{xx}v_n\Sigma \Sigma^\top] + \partial_xv_n\Sigma\cdot q  \right)\diff t \\
		&\quad + \int_s^1\partial_xv_n\Sigma\big(t, X(t)\big)\cdot\diff W(t)\\
		&= \int_s^1\partial_xv_n\Sigma\cdot q - g_n\big(t,q(t,X(t)),X(t) \big)\diff t + \int_s^1\partial_xv_n\Sigma\big(t, X(t) \big)\cdot\diff W(t),
	\end{align*}
	with the functions in the drifts being evaluated at $(t,X(t))$, and where the second equality follows by \eqref{eq:PDE.bsde.n} and \eqref{eq:subgrad}.
	Taking conditional expectations on both sides and noticing that by Fenchel-Moreau theorem $g_n\ge g$, one has
	\begin{align}\notag
		\rho^{g_n}_s\big(F(X^{s,\xi}(1)\big)& = v_n(s, \xi ) 
		 = \E\bigg[F\big( X(1) \big) - \int_s^1g_n\big(t,q(t,X(t)), X(t)\big)\,dt \bigg|\cF_s\bigg]\\
		 \label{eq:upper.bound.n.prop}
		  &\le \esssup_{q  \in \cL^\infty}\E\bigg[F(X^{s,\xi,q}(1)) - \int_s^1g\big(t,q(t),X^{s,\xi,q}(t)\big)\diff t\bigg|\cF_s \bigg].
	\end{align}
	Since $g_n\uparrow g$ pointwise and $\lim_n\rho_0^{g_n}\big( F(X^{s,x}(1)) \big)\le K$, (this follows from the fact that $\|F\|\le K$ and $g\ge0$) it follows from \cite[Theorem 4.14]{DHK1101} that
	\begin{equation*}
		\rho^{g_n}_s\big( F(X^{s,\xi}(1) \big)\uparrow \rho^g_s\big( F(X^{s,\xi}(1) \big) \quad\P\text{--a.s.}
	 \end{equation*} 
	Hence, taking the limit as $n$ goes to infinity, it holds
	\begin{align*}\
		&\rho_s^{g}\big(F(X^{s,x}(1)\big) \le \esssup_{ q\in {\cal L}^\infty}\E\bigg[F\big(X^{s,x,q}(1) \big) - \int_s^1g\big( t,q(t),X^{s,x,q}(t) \big)\diff t \bigg|\cF_s \bigg].
	\end{align*}
	When $F$ is not bounded, apply the result to $F^m:= F\wedge m$ and derive the above upper bound for $F$ by monotone convergence.
	Combined with Step 1, this yields the result.
\end{proof}

\begin{remark}
\label{rem:Linear.growth.case}
	It is interesting to notice that, from the proof of Proposition \ref{prop:var.rep}, when the function $f$ is of linear growth, the supremum in \eqref{eq:var.rep.abstract.} can be restricted to $q$ bounded by a fixed constant.
	In fact, take for instance the truncated function $f_n$.
	Then, by \eqref{eq:upper.bound.lips.n} and \eqref{eq:upper.bound.n.prop} we have
	\begin{align*}\
		\rho_s^{g_n}\big(F(X^{s,x}(1)\big) = \sup_{ q\in \cL^\infty_n}\E\bigg[F\big(X^{s,x,q}(1) \big) - \int_s^1g_n\big( t,q(t),X^{s,x,q}(t) \big)\diff t \bigg|\cF_s \bigg]
	\end{align*}
	where $\cL^\infty_n$ is the set of elements of $\cL^\infty$ bounded by $n$.
\end{remark}
The proof of Theorem \ref{thm:var.rep.McKV} will need the ensuing lemma.
Essentially, this lemma is only a slightly modified version of \cite[Lemma A.1]{BaLaTa}.
We state it separately since it will be used repeatedly in the subsequent sections. 

\begin{lemma}[Tightness lemma]
\label{Lem:tightness}
	Let the conditions of Theorem \ref{thm:var.rep.McKV} be satisfied.
	Let $(q_n)_{n\ge1}$ be a sequence in $\cL^\infty$ such that
	\begin{equation*}
		\E\bigg[\int_0^1g_1(u, q_n(u))\diff u\bigg]\le C
	\end{equation*}
	for some constant $C>0$ and let $\beta$ be a bounded, predictable process with values in $\R^{\xdim\times d}$.
	Then $(q_n)_{n\ge1}$ admits a subsequence converging to $q$ in the weak topology of $\L^2([0,T]\times\Omega)$.
	Putting $A_n:= \int_0^\cdot \beta (u)q_n(u)\diff u$,
	it holds that the sequence $(A_n)_{n\ge1}$ is tight, there is a continuous process $A$ and a subsequence $(A_{n_k})_{k\ge1}$ such that $(A_{n_k})_{k\ge1}$ converges to $A$ in law in $\cC([0,T],\R^\xdim)$.
	Furthermore, it holds $A(t) = \int_0^t\beta(u)q(u)\diff u$ and
	\begin{equation*}
		\liminf_{n\to \infty}\E\bigg[ \int_0^1g_1(u, q_n(u))\diff u\bigg]\ge \E\bigg[ \int_0^1g_1(u, q(u))\diff u\bigg].
	\end{equation*}
\end{lemma}
\begin{proof}
	The case where $\beta$ is a constant is exactly \cite[Lemma A.1]{BaLaTa}.
	In particular, this lemma shows that a subsequence of $(q_n)_{n\ge1}$ admits a weak limit in $\L^2([0,T]\times \Omega)$.
	The tightness of $A_n(t)=\int_0^t\beta(u)q_n(u)\diff u$ follows exactly as in the proof of \cite[Lemma A.1]{BaLaTa} since $\beta$ is bounded.
	Therefore, $(A_n)_{n\ge1}$ admits a subsequence $(A_{n_k})_{k\ge 1}$ converging to some $A$ in law, and therefore $\P$-a.s. on some probability space $(\Omega,\cF, \P)$.
	It remains to show that $A_t = \int_0^t\beta(u)q(u)\diff u$.
	First, notice that by the coercivity condition $(A1)$, the convex conjugate $g_1$ of $f_1$ also satisfies $\lim_{\|q\|\to \infty}\inf_{t \in [0,1]}\frac{g_1(t,q)}{\|q\|} = \infty$, see e.g. \cite[Theorem 3.3]{Evans1998}.
	Thus, there is $\varepsilon>1$ and $a>0$ such that for $\|q\|$ large enough, $g_1(t,q)\ge a\|q\|^\varepsilon$.
	Thus,
	\begin{equation*}
		\E\bigg[\int_0^1\|\beta(u)q_n(u)\|^\varepsilon\diff u\bigg] \le (C_1\|\beta\|_\infty)^\varepsilon+ \frac{\|\beta\|^\varepsilon_\infty}{a}\E\bigg[\int_0^1g_1(u, q_n(u))\diff u\bigg] \le C_2
	\end{equation*}
	for some $C_1, C_2>0$.
	Thus, by the de la Vall\'ee Poussin compactness criterion, the sequence $(\beta q_n)_{n\ge1}$ is weakly relatively compact in $\L^1([0,1]\times \Omega)$.
	Thus, up to a subsequence, it converges to some $\beta q$.
	To show that $A_t = \int_0^t\beta(u)q(u)\diff u$, let $Z\in \L^\infty(\R^\xdim,\cF_T)$.
	Thus,
	\begin{equation*}
		\E[ZA(t)] = \lim_{n\to \infty}\E[ZA_n(t)] = \lim_{n\to \infty}\E\bigg[Z\int_0^t\beta(u)q_n(u)\diff u\bigg] = \E\bigg[Z\int_0^t\beta(u)q(u)\diff u\bigg],
	\end{equation*}
	which shows that $A_t=\int_0^t\beta(u)q(u)\diff u$ for all $t\in [0,1]$ and by continuity of both processes, $A= \int_0^\cdot \beta(u)q(u)\diff u$.
\end{proof}
We conclude this section with the proof of Theorem \ref{thm:var.rep.McKV}, the variational representation for functions of McKean-Vlasov equations.
\begin{proof}[Proof of Theorem \ref{thm:var.rep.McKV}]
	We first assume that $F\in C_b(\R^\xdim,\cP_2(\R^\xdim))$, and we consider the particle system
	\begin{align}
	\label{eq:N.McKv}
		X^{i,N}(t) = x + \int_s^tb\big(u, X^{i,N}(u), L^N(\X(u)) \big)\diff u + \int_0^t\sigma\big(u, X^{i,N}(u), L^N(\X(u)) \big)\diff W^i(u)
	\end{align}
	for $N$ $d$-dimensional independent Brownian motions $\W:=(W^1,\dots, W^N)$.
	We assume that $W^1=W$, the driving Brownian motion in \eqref{eq:contr.MckV}.
	Let us introduce the functions
	\begin{align}
	\label{eq:def.B.Sigma}
		B(t,\x) &:= \Big( b(t, x^i, L^N(\x)) \Big)_{i=1,\dots,N},\quad \Sigma(t,\x) := \text{diag}\bigg(\Big( \sigma(t, x^i, L^N(\x)) \Big)_{i=1,\dots,N}\bigg)
	\end{align}
	\begin{align*}
		\widehat g(t,q,\x) &:= \Big( g(t,q,x^i,L^N(\x))\Big)_{i=1,\dots,N},\quad \widehat F(\x) := \Big( F(x^i,L^N(\x)) \Big)_{i=1,\dots,N},\quad (t, \x) \in [0,1]\times (\R^\xdim)^N.
	\end{align*}
	Then, the vector $\X:= (X^{1,N},\dots, X^{N,N})$ satisfies the equation
	\begin{equation*}
		\X(t) = \boldsymbol{\x} + \int_s^tB(u, \X(u))\diff u + \int_s^t\Sigma(u, \X(u))\diff \W(u)
	\end{equation*}
	and $\x = (x,\dots,x)\in (\R^\xdim)^N$.
	Since $b$ and $\sigma$ are Lipschitz--continuous in the measure argument with respect to the Wasserstein distance, it follows that $B$ and $\Sigma$ are Lipschitz--continuous with respect to the Euclidean norm, with the same Lipschitz constant $\ell_b$, and Equation \eqref{eq:N.McKv} is thus well-posed.
	Let us assume for the moment that $b$ is bounded and observe that by the conditions on $f$ and $F$, the functions $\widehat F$ and $\widehat g$ are continuous.
	Thus, the functions $B$, $\Sigma$, $\widehat F$ and $\widehat g$  satisfy the conditions of Proposition \ref{prop:var.rep} from which we get
	\begin{align}
	\notag
		&\rho^{g}_s\big(F(X^{i,N}(1), L^N(\X(1)) \big) = \rho^{\widehat g}_s(\widehat F(\X^{s,\xi}(1)))\\\notag
		&\qquad = \esssup_{q \in \cL^\infty}\E\bigg[\widehat F\big(\X^{s,x,q}(1) \big) - \int_s^1\widehat g\big(s, q(u), \X^{s,x,q}(u)\big)\diff u \bigg|\cF_s \bigg]\\
		\label{eq:rep.N}
		& \qquad=\esssup_{q \in \cL^\infty}\E\bigg[ F\big(\X^{i,N,q}(1), L^N(\X^{q}(1)) \big) - \int_s^1g\big(u, q(u), X^{i,N,q}(u), L^N(\X^{q}(u))\big)\diff u \bigg|\cF_s \bigg]
	\end{align}
	where in the second line above we omitted the superscript $(s,\xi)$ to simplify the notation.
	The goal is now to take the limit on both sides as $N$ goes to infinity, using propagation of chaos arguments.
	However, the convergence of the left hand side is hard in view of the weak growth condition assumed on $g_1$.
	To overcome this issue, we will again truncate this function, and first prove the limit for the truncated function, which is Lipschitz--continuous.
	Thus, we further consider the functions $f_n$ and $g_n$ defined in \eqref{eq:trunc.f} and \eqref{eq:conjug.trunc} respectively.
	By Proposition \ref{prop:var.rep} (or actually Remark \ref{rem:Linear.growth.case}),
	\begin{align}
	\notag
		&\rho^{g_n}_s\big(F(X^{i,N}(1), L^N(\X(1)) \big)\\
		\label{eq:rep.N.n}
		&  \quad=\esssup_{q \in \cL_n^\infty}\E\bigg[ F\big(\X^{i,N,q}(1), L^N(\X^{q}(1)) \big) - \int_s^1g_n\big(u, q(u), X^{i,N,q}(u), L^N(\X^{q}(u))\big)\diff u \bigg|\cF_s \bigg].
	\end{align}

	Let $n\ge1$ be fixed.
	The function $f_n$ is Lipschitz--continuous in its second variable and as argued in the proof of Proposition \ref{prop:var.rep}, $\rho^{g_n}_s\big(F(X^{i,N}(1), L^N(\X(1))) \big) = Y_n^N(s)$ where $(Y_n^N, Z^N_n)$ solves the BSDE \eqref{eq:bsde} with the terminal condition $F(X^{i,N}(1), L^N(\X(1)))$ and the generator $(s,z)\mapsto f_n\big(s, z, X^{i,N}(s), L^N(\X(s))\big)$.
	By standard propagation of chaos results, see e.g. \cite[Theorem 2.12]{MR3753660}, the sequence $(X^{1,N})_{N\ge 1}$ converges to $X^{s,\xi}$ in $\mathbb{S}^2(\R^\xdim,\F)$ and the sequence $(L^N(\X(t)))_{N\ge1}$ converges to $\mu^{s,\xi}(t)$ in second order Wasserstein distance.
	By continuity of the functions $F$ and $f_n$,
	it follows that $F(X^{i,N}(1), L^N(\X(1)))$ converges to $F(X^{s,\xi}(1), \mu^{s,\xi}(1))$ in $\L^2(\R^\xdim,\cF_1)$ and $f_n\big(s, z, X^{i,N}(s), L^N(\X(s))\big)$ converges to $f_n\big(s, z, X^{s,\xi}(s), \mu^{s,\xi}(s)\big)$ in $\H^2(\R,\cF)$.
	Thus, using stability for BSDE solutions with Lipschitz--continuous generators, see \cite[Proposition 2.1]{karoui01}, we obtain that (up to a subsequence) $(Y^N_n(s))_{N\ge1}$ converges to  $\rho^{g_n}_s\big(F(X^{s,\xi}(1), \mu^{s,\xi}(1))) \big) = Y_n(s)$, where $(Y_n,Z_n)$ is the solution of the BSDE \eqref{eq:bsde} with terminal condition $F(X^{s,\xi}(1), \mu^{s,\xi}(1))$ and generator $(s,z)\mapsto f_n\big(s, z, X^{s,\xi}(s), \mu^{s,\xi}(s)\big)$.
	Taking the limit in \eqref{eq:rep.N.n} thus yields, for all $q\in {\cal L}^\infty_n$
	\begin{align}
	\label{eq:limit.rho.N}
	&\rho^{g_n}_s\big(F(X^{s,\xi}(1), \mu^{s,\xi}(1)) \big)=
		\lim_{N\to \infty}\rho^{g_n}_s\big(F(X^{i,N}(1), L^N(\X(1))) \big)\\
		\notag
		&\qquad\qquad \ge \liminf_{N\to \infty}\E\bigg[F\big(\X^{i,N,q}(1), L^N(\X^{q}(1)) \big) - \int_s^1g_n\big(s, q, X^{i,N,q}(u), L^N(\X^{q}(u))\big)\diff u \bigg|\cF_s \bigg]\\
		\label{eq:McKV.rep.n.lowerbound}
		&\qquad\qquad \ge \E\bigg[F\big(\X^{s,\xi,q}(1), \mu^{s,\xi,q}(1)) \big) - \int_s^1g_n\big(s, q(s), X^{s,\xi,q}(u), \mu^{s,\xi,q}(u)\big)\diff u \bigg|\cF_s \bigg]
	\end{align}
	where the second inequality follows by dominated convergence and propagation of chaos, since the interacting particle system
	\begin{align*}
	\label{eq:M,q.McKv}
		X^{i,N,q}(t) &= \xi + \int_s^tb\big(u, X^{i,N,q}(u), L^N(\X^q(u)) \big) + \sigma\big(u, X^{i,N,q}(u), L^N(\X^q(u)) \big)q(u)\diff u\\
		 &\quad+ \int_0^t\sigma\big(u, X^{i,N,q}(u), L^N(\X^q(u)) \big)\diff W^i(u)
	\end{align*}
	converges (in the sense explained above for $X^{i,N}$) to the McKean-Vlasov equation \eqref{eq:contr.MckV}, see \cite[Theorem 1.12]{MR3753660}.
	Since $n\ge1$ and $q \in \cL^\infty_n$ were taken arbitrary, letting $n$ go to infinity, using the arguments leading to \eqref{eq:upper.bound.lips}, and then taking the supremum over $q \in \cL^\infty$, we obtain
	\begin{equation}
	\label{eq:MckV.rep.lower}
	\rho^{g}_s\big(F(X^{s,\xi}(1), \mu^{s,\xi}(1)) \big)
		\ge \esssup_{q\in \cL^\infty}\E\bigg[F\big(\X^{s,\xi,q}(1), \mu^{s,\xi,q}(1)) \big) - \int_s^1g\big(s, q, X^{s,\xi,q}(u), \mu^{s,\xi,q}(u)\big)\diff u \bigg|\cF_s \bigg].
	\end{equation}

\vspace{.2cm}

	Let us now derive the upper bound.
	Let again $n$ be fixed.
	By \eqref{eq:rep.N.n} for every $N \ge 1$, there is $q_N\in \cL^\infty_n$ such that
	\begin{align}
	\notag
		&\rho^{g_n}_s\big(F(X^{i,N}(1), L^N(\X(1)) \big) \\
		\label{eq:ineq.n.N.upper}
		&\le \E\bigg[ F\big(X^{i,N,q_N}(1), L^N(\X^{q_N}(1)) \big) - \int_s^1g_n\big(u, q_N(u), X^{i,N,q_N}(u), L^N(\X^{q_N}(u))\big)\diff u\bigg| \cF_s\bigg] + \frac1N.
	\end{align}
	The construction of such an approximate optimizer is classical, it is detailed for instance in the proof of Corollary \ref{cor:F.W-random} below.
	Obviously, $q_N$ depends on $n$, but this will not play any role in what follows.
	By boundedness of $F$ and $g_2$ there is a constant $C>0$ such that
	\begin{equation*}
		\E\bigg[\int_0^1g_{1}(u,q_N(u))\diff u \bigg]\le C\quad \text{for all } N\ge1.
	\end{equation*}
	Hence, by Lemma \ref{Lem:tightness} and the boundedness of $\sigma$, there is some $q\in \cL$ such that the sequences of processes $\int_0^\cdot q_N(u)\diff u$ and $\int_0^\cdot q_N(u)\sigma(u, X^{s,\xi,q}(u),\mu^{s,\xi,q}(u))\diff u$ converge in law to the processes $\int_0^\cdot q(u)\diff u$ and $\int_0^\cdot q(u)\sigma(u, X^{s,\xi,q}(u),\mu^{s,\xi,q}(u))\diff u$ respectively,	where $X^{s,\xi,q}$ is the solution of the McKean-Vlasov equation \eqref{eq:contr.MckV} and
	\begin{equation}
	\label{eq:liminf.N.chaos}
		\liminf_{N\to \infty}\E\bigg[ \int_0^1g_1(u, q_N(u))\diff u\bigg]\ge \E\bigg[ \int_0^1g_1(u, q(u))\diff u\bigg].
	\end{equation}
	Let us now show that $(X^{i,N,q_N}_t)_{N\ge 1}$ converges to $X^q_t$ in $\L^2(\R^\xdim, \cF_t)$, possibly on another probability space.
	This follows again by propagation of chaos arguments which we give for completeness.
	Since the sequence $\int_0^\cdot q_N(u)\sigma(u, X^{s,\xi,q}(u),\mu^{s,\xi,q}(u))\diff u$ converges in law in $\cC([0,T],\R^\xdim)$ to $\int_0^\cdot q(u)\sigma(u, X^{s,\xi,q}(u),\mu^{s,\xi,q}(u))\diff u$, it follows by Skorohod's representation theorem that, passing to a subsequence again indexed by $N$, we have
	\begin{equation*}
		\int_0^\cdot q_N(u)\sigma(u, X^{s,\xi,q}(u),\mu^{s,\xi,q}(u))\diff u\longrightarrow \int_0^\cdot q(u)\sigma(u, X^{s,\xi,q}(u),\mu^{s,\xi,q}(u))\diff u \quad \P\text{--a.s.}
	\end{equation*}
	with all processes defined on a common probability space again denoted $(\Omega, \cF, \P)$.
	Let $(\widetilde X^{i,q})_{i\ge1}$ be i.i.d. copies of $X^{s,\xi,q}$ such that $\widetilde X^{i,q}$ satisfies \eqref{eq:contr.MckV} with the driving Brownian motion $W^i$.
	In particular, by the law of large numbers, we have that $\cW^2_2\big(L^N(\widetilde \X^q(t)), \mu^q(t)\big)$ converges to zero $\P$--a.s., uniformly in $t$.
	Then, by Lipschitz--continuity of $b$ and $\sigma$ and using triangular inequality and Gronwall's inequality,
	\begin{align}
	\notag
		&\E\Big[\| X^{i,N,q_N}(t) - \widetilde X^{s,\xi,q}(t)\|^2 \Big]\\\notag
		& \le e^{L_b^2(2 + \|q_N\|_\infty^2)}L_b(2 + \|q_N\|_\infty)\bigg\{ \E\bigg[\int_0^t\frac1N\sum_{j=1}^N\|X^{j,N,q}(u) - \widetilde X^{s,\xi,q}(u)\|^2\diff u\bigg]\\
		\label{eq:chaos.gronwall}
		&\quad + \E\bigg[\bigg\|\int_0^1\sigma\big(u,X^{s,\xi,q}(u), \mu^{s,\xi,q}(u)\big)\big\{ q_N(u) - q(u)\big\}\diff u\bigg\|^2 + \int_0^1 \cW^2_2\big(L^N(\widetilde \X^{s,\xi,q}(u), \mu^{s,\xi,q}(u)\big)\diff u  \bigg] \bigg\}.
	\end{align}
	Since $\|q_N\|_\infty\le n$,
	averaging on both sides and applying again Gronwall's inequality,
	\begin{align*}
		\frac1N\sum_{i=1}^N\E\Big[\| X^{i,N,q_N}(t) - \widetilde X^{s,\xi,q}(t) \Big]
		&\le C_{L_b,n} \E\bigg[\bigg\|\int_0^1\sigma(u,X^{s,\xi,q}(u), \mu^{s,\xi,q}(u))\big\{q_N(u) - q(u)\big\}\diff u\bigg\|^2\\
		&\qquad + \int_0^1 \cW^2_2\big(L^N(\widetilde \X^{s,\xi,q}(u)), \mu^q(u)\big)\diff u  \bigg]\xrightarrow[N\to\infty]{} 0.
	\end{align*}
	Plugging this back into \eqref{eq:chaos.gronwall} thus implies by dominated convergence that $(X^{1,N,q_N}(t))_{N\ge1}$ converges to $X^{s,\xi,q}(t)$ in $\L^2(\R^\xdim,\cF)$.
	Applying triangular inequality,
	\begin{equation*}
		\E\Big[ \cW^2_2\big( L^N(\X^{q_N}(t)), \mu^q(t) \big) \Big] \le 2 \E\bigg[ \cW^2_2\big( L^N(\widetilde\X^q(t)), \mu^q(t) \big) + \frac1N\sum_{j=1}^N\|X^{j,N,q_N}(t) - \widetilde X^{j,N,q}(t)\|^2 \bigg].
	\end{equation*}
	Thus, $\big( L^N(\X^{q_N}(t)) \big)_{N\ge 1}$ converges to $\mu^q(t)$ in second order Wasserstein distance.
	Coming back to \eqref{eq:ineq.n.N.upper}, taking the limit $N\longrightarrow \infty$ therein and recalling that $g_n\ge g$, it follows by \eqref{eq:liminf.N.chaos} and \eqref{eq:limit.rho.N} that
	\begin{align}
	\notag
		&\rho^{g_n}_s\big(F(X^{s,\xi}(1), \mu^{s,\xi}(1)) \big) \le \E\bigg[ F\big(X^{q}(1),\mu^{s,\xi,q}(1) \big) - \int_s^1g\big(u, q(u), X^{s,\xi,q}(u), \mu^{s,\xi,q}(u)\big)\diff u \bigg|\cF_s\bigg]\\
		\label{eq:McKV.rep.n}
		&\le \esssup_{q \in \cL}\E\bigg[ F\big(X^{s,\xi,q}(1),\mu^{q}(1) \big) - \int_s^1g_n\big(u, q(u), X^{s,\xi,q}(u), \mu^{s,\xi,q}(u)\big)\diff u\bigg|\cF_s\bigg]\\\notag
		&\le \esssup_{q \in \cL}\E\bigg[ F\big(X^{s,\xi,q}(1),\mu^{q}(1) \big) - \int_s^1g\big(u, q(u), X^{s,\xi,q}(u), \mu^{s,\xi,q}(u)\big)\diff u\bigg|\cF_s\bigg].
	\end{align}
	Since $g_n\uparrow g$ pointwise and $\lim_n\rho_0^{g_n}\big( F(X^{s,x}(1), \mu^{s,x}(1))) \big)\le K$, it holds 
	\begin{equation*}
		\rho^{g_n}_s\big( F(X^{s,\xi}(1) \big)\uparrow \rho^g_s\big( F(X^{s,\xi}(1) \big) \quad\P\text{--a.s.},
	 \end{equation*} 
	 see \cite[Theorem 4.14]{DHK1101}.
	 Hence,
	 \begin{equation*}
	 	\rho^{g}_s\big(F(X^{s,\xi}(1), \mu^{s,\xi}(1)) \big)
		\le \esssup_{q \in \cL}\E\bigg[ F\big(X^{s,\xi,q}(1),\mu^{s,\xi,q}(1) \big) - \int_s^1g\big(u, q(u), X^{s,\xi,q}(u), \mu^{q}(u)\big)\diff u\bigg|\cF_s\bigg].
	 \end{equation*}
	Combine this with \eqref{eq:MckV.rep.lower} to get \eqref{eq:var.rep.abstract.McKV} and conclude the proof for the case $b$ bounded.

	When $b$ is unbounded, we approximate it by a sequence of bounded, Lipschitz--continuous function $b_k$ and, denoting by $X_k^{s,\xi,q}$ the solution of the McKean--Vlasov SDE \eqref{eq:contr.MckV} with drift $b$ replaced by $b_k$, it follows by \eqref{eq:McKV.rep.n.lowerbound} and \eqref{eq:McKV.rep.n} that
	\begin{align}
	\notag
	 	&\rho^{g_n}_s\big(F(X_k^{s,\xi}(1), \mu_k^{s,\xi}(1)) \big)\\\label{eq:rep.n.k.limit}
		&= \esssup_{q \in \cL}\E\bigg[ F\big(X_k^{s,\xi,q}(1),\mu_k^{s,\xi,q}(1) \big) - \int_s^1g_n\big(u, q(u), X_k^{s,\xi,q}(u), \mu_k^{s,\xi,q}(u)\big)\diff u\bigg|\cF_s\bigg].
	\end{align}
	Now, observe that as $k$ goes to infinity, it follows by stability of McKean-Vlasov SDEs that the sequences $(X_k^{s,\xi,q})_{k\ge1}$ and $(\mu_k^{s,\xi,q})_{k\ge1}$ converge, respectively, to $X^{s,\xi,q}$ and $\mu^{s,\xi,q}$, (solution of \eqref{eq:contr.MckV}) in $\S^2(\R^\xdim,\F)$ and in Wasserstein distance, respectively.
	Therefore, using exactly the same arguments developed above allows to take the limit as $n$ and $k$ go to infinity in \eqref{eq:rep.n.k.limit} to obtain the desired representation \eqref{eq:var.rep.abstract.McKV}.

	When $F$ is only bounded from below and lower semicontinuous (and since $\R^\xdim\times \cP_2(\R^\xdim)$ is a metric space) it can be approximated by an increasing sequence of bounded continuous functions $F^n$.
	We can thus apply the representation \eqref{eq:var.rep.abstract.McKV} to $F^n$ and obtain the result for $F$ by monotone convergence.
\end{proof}

\section{Scaling limits for functionals of McKean-Vlasov diffusions}
The representation theorem given in Theorem \ref{thm:var.rep.McKV} will allow us to derive two limit theorems for the functional $\rho^g(F(X^{s,\xi}(1), \mu^{s,\xi}(1)))$ when the noise coefficient $\sigma$ is scaled to zero.
These results are essentially generalizations of Freidlin-Wentzell theorem in large deviations theory to the non-exponential case, while while considering McKean-Vlasov diffusion. 
Let us recall that a similar generalization of Schilder's theorem is given in \cite{BaLaTa}.

\subsection{A non-exponential Freidlin-Wentzell theorem}
The goal of this subsection is to prove the following result.
We will argue at the end of the section that Corollary \ref{cor:F-WThm} follows as a direct consequence.
Recall that the space $\mathcal{H}$ is defined in \eqref{eq:def.calH}.
 
\begin{theorem}
\label{thm:abstract.schilder}
	Assume that the conditions $(A1)$, $(A2)$ and $(A3)'$ are satisfied.
	Put 
	\begin{equation}
	\label{eq:def.g.n}
		g_n(t,q, x,\mu):= g\big(t, \frac{q}{\sqrt{n}}, x, \mu\big)
	\end{equation}
	and given $s \in [0,1]$ and $x \in \R^\xdim$, let $X_n$ solve the SDE
	\begin{equation}
	\label{eq:McKV.SDE.n}
		\begin{cases}
		dX_n(t) = b\big(t, X_n(t), \mu_n(t)\big)\diff t + \frac{1}{\sqrt{n}}\sigma\big(t, X_n(t), \mu_n(t)\big)\diff W(t),\\
		 \quad X_n(s) = x,\,\, \mu_n(t) = \text{law}(X_n(t)).
		\end{cases}
	\end{equation}
	For every $F \in C_b(\mathbb{R}^m\times \cP_2(\mathbb{R}^m))$ it holds\footnote{By convention, we set $\int_s^1g\big(t, \varphi(t), \Phi(t),\delta_{\Phi(t)} \big)\diff t = +\infty$ when $\varphi$ is such that \eqref{eq:ODE} does not have a solution.}
	\begin{equation}
	\label{eq:abstract.schilder}
		\lim_{n\to \infty}\rho^{g_n}_s\big(F(X_n(1), \mu_n(1))\big) = \sup_{ \varphi \in  \mathcal{H}}\bigg(F\big( \Phi^\varphi(1),\delta_{\Phi^\varphi(1)} \big) - \int_s^1g\big(t, \varphi(t), \Phi^\varphi(t),\delta_{\Phi^\varphi(t)} \big)\diff t  \bigg)
	\end{equation}
	where $\Phi^\varphi$ is the solution (when it exists) of the ordinary differential equation
	\begin{equation}
	\label{eq:ODE}
			\diff \Phi^\varphi(t) =  b\big(t, \Phi^\varphi(t), \delta_{\Phi^\varphi(t)}\big) + \sigma\big(t, \Phi^\varphi(t), \delta_{\Phi^\varphi(t)} \big)\varphi(t)\diff t,\quad \Phi^\varphi(s) = x. 
	\end{equation}
\end{theorem}

\begin{proof}
	The starting point of the proof is the general variational representation \eqref{eq:var.rep.abstract.McKV} which, in the case of $g_n$ and $X_n$ reads
	\begin{align}
	\notag
		\rho_s^{g_n}\big(F(X_n(1),\mu_n(1)) \big)
		&= \esssup_{q\in \mathcal{L} }\E\bigg[ F\big(X^{\frac{q}{\sqrt{n}}}_n(1), \mu^{\frac{q}{\sqrt{n}}}_n(1)\big) - \int_s^1g_n\big(u,q_u,X^{\frac{q}{\sqrt{n}}}_n(u),\mu^{\frac{q}{\sqrt{n}}}_n(u))\big)\diff u \bigg| \cF_s\bigg]\\
		\label{eq:rep.proof.schilder}
		& = \esssup_{q\in \mathcal{L} }\E\bigg[F\big(X^{q}_n(1), \mu^{q}_n(1)\big) - \int_s^1g\big(u,q(u),X^{q}_n(u),\mu^{q}_n(u))\big) \diff u \bigg|\cF_s \bigg]
	\end{align}
	where we omitted the superscript $(s,\xi)$ in $X_n$ to make the notation lighter.
	We will show the upper and lower bounds in \eqref{eq:abstract.schilder} separately:

\vspace{.2cm}

	\emph{Step 1: The lower bound}.
	We first show that for every $\varphi \in \mathcal{H}$ it holds $X^\varphi_n(t) \rightarrow \Phi(t)$ $\P$--a.s. and $\cW_2(\mu^\varphi_n(t) , \delta_{\Phi(t)})\rightarrow 0$.
	By the growth property of $b$ and boundedness of $\sigma$, we have
	\begin{align*}
		\E\Big[\|X^\varphi_n(t)\|^2 \Big] \le 8\|x\|^2 + 32\ell_b^2(1 + \frac{4}{n}) +8\ell^2_b\|\varphi\|^2_1 + 64\ell_b^2 \int_0^t\E\Big[\|X^\varphi_n(u)\|^2 \Big]\diff u,
	\end{align*}
	where we also used It\^o's isometry, and where we used the notation 
	$$\|\varphi\|_p:= \bigg(\int_0^1\|\varphi(u)\|^p\diff u\bigg)^{1/p}, \quad p\ge1.$$
	Thus, it follows by Gronwall's inequality that
	\begin{equation}
	\label{eq:bounded.moment}
		\E\Big[\|X^\varphi_n(t)\|^2 \Big] \le e^{64\ell_b^2}\Big(8\|x\|^2 + 32\ell_b^2(1 +\frac4n) +8\ell^2_b\|\varphi\|^2_1\Big).
	\end{equation}
	That is, the sequence $(X_n^\varphi)_{n\ge1}$ is bounded in $\L^2(\R^\xdim, \cF_t)$.
	Now by Lipschitz--continuity of $b$ and $\sigma$ and Burkholder-Davis-Gundy's inequality, it follows that 
	\begin{align}
	\notag
		\E\Big[\|X^\varphi_n(t) - \Phi^\varphi(t)\|^2 \Big] &= 4 (\ell_b^2 + \|\varphi\|_2^2)\E\bigg[\int_0^t\|X_n^\varphi(u) - \Phi(u)\|^2 + \cW_2^2(\mu^\varphi_n(u), \delta_{\Phi^\varphi(u)})\diff u \bigg]\\
		\label{eq:boud.varphi.detem}
		&\quad + \frac{\ell_b^2}{n}\E\bigg[\int_0^t1 + \|X^\varphi_n(u)\|^2 + \E\big[\|X^\varphi_n(u)\|^2\big]\diff u \bigg].
	\end{align}
	By \eqref{eq:bounded.moment}, there is a constant $C>0$ such that the last term is bounded by $C/n$ and by properties of the Wasserstein distance, it holds that $\cW_2^2(\mu^\varphi_n(u), \delta_{\Phi^\varphi(u)}) \le \E[\|X^\varphi_n(u) - \Phi^\varphi(u)\|^2]$.
	Therefore, applying Gronwall's inequality gives
	\begin{align*}
		\E\Big[\|X^\varphi_n(t) - \Phi^\varphi(t)\|^2\Big] \le e^{4(\ell_b^2 + \|\varphi\|_2^2)}\frac{C\ell_b^2}{n},
	\end{align*}
	which shows $\P$-a.s. convergence of $X^\varphi_n(t)$ to $\Phi^\varphi(t)$, at least for a subsequence, and also convergence of $\cW_2^2(\mu^\varphi_n(u), \delta_{\Phi^\varphi(u)})$ to zero.

	To prove the lower bound, let $ \varphi$ be an arbitrary element of $\mathcal{H}$, then $\varphi \in \cL$.
	By \eqref{eq:rep.proof.schilder}, it follows by dominated convergence and continuity of $F$ and $g$ in its last two components that
	\begin{align*}
		\liminf_{n\to \infty}\rho^{g_n}_s\big( F(X_n(1), \mu_n(1)) \big)
		& \ge \liminf_{n\to \infty}\E\bigg[ F\big(X^{\varphi}_n(1), \mu^{\varphi}_n(1)\big) - \int_s^1g\big(u,\varphi(u),X^{\varphi}_n(u),\mu^{\varphi}_n(u)\big)\diff u \bigg| \cF_s\bigg]\\
		&\ge F(\Phi^\varphi(1), \delta_{\Phi^\varphi(1)}) - \int_s^1g\big(u, \varphi(u), \Phi^\varphi(u), \delta_{\Phi^\varphi(u)}\big)\diff u .
	\end{align*}
	Since $\varphi$ was taken arbitrary, this proves the lower bound in \eqref{eq:abstract.schilder}.

\vspace{.2cm}

	\emph{Step 2: The upper bound}.
	Since $X_n(s)=x$ is deterministic, it follows that $\rho_s^{g_n}\big(F(X_n(1), \mu_n(1) \big)$ is deterministic as well.
	Thus, \eqref{eq:rep.proof.schilder} becomes
\begin{align}
	\notag
		\rho_s^{g_n}\big(F(X_n(1),\mu_n(1)) \big)
		& = \sup_{q\in \mathcal{L} }\E\bigg[F\big(X^{q}_n(1), \mu^{q}_n(1)\big) - \int_s^1g\big(u,q(u),X^{q}_n(u),\mu^{q}_n(u) \big) \diff u \bigg].
	\end{align}	
	Let $n \in \mathbb{N}$ and let $ q_n$ come $1/n$--close to the optimal value above, i.e. 
	\begin{align}
	\label{eq:upper.bound.n.close}
		\rho_s^{g_n}\big(F(X_n(1),\mu_n(1)) \big)
		 \le \E\bigg[ F\big(X^{q_n}_n(1), \mu^{q_n}_n(1))\big) - \int_s^1g\big(u,q_n(u),X^{q_n}_n(u),\mu^{q_n}_n(u)\big)\diff u \bigg] + \frac1n.
	\end{align}
	In particular,
	\begin{equation}
	\label{eq:bound.g.lower}
		\E\bigg[ \int_s^1g_1(u, q_n(u))\diff u \bigg]\le C\quad \text{for all } n\ge1\quad \text{for some } C>0.
	\end{equation}
	Hence, due to the coersivity condition in assumption $(A1)$, it follows by Lemma \ref{Lem:tightness} that (up to a subsequence), the sequence $(q_n)_{n\ge1}$ converges weakly to some $q\in \cL$, the sequence $A_n = \int_0^\cdot q_n(u)\diff u$ is tight in $C([0,1],\R^d)$ and it admits a subsequence converging to $A = \int_0^\cdot q(u)\diff u$.
	Moreover,
	\begin{equation}
	\label{eq:fatou.upper}
		\liminf_{n\to \infty}\E\bigg[\int_s^1g_1(t, q_n(t))\diff t \bigg] \ge \E\bigg[\int_s^1g_1(t,  q(t))\diff t \bigg].
	\end{equation}
	Next, we will show that the sequence  of processes $(X^{q_n}_n)$ is tight and admits a subsequence converging to $X^q$ in law on $C([0,T],\R^\xdim)$.
	
	Let us first show the tightness property.
	For all $0\le s\le t\le 1$, it follows by linear growth of $b$ and boundedness of $\sigma$ that
	\begin{align}
	\notag
		\E\Big[\|X^{q_n}_n(t) - X^{q_n}_n(s)\|\Big] &\le L_b\int_s^t\Big(1 + \E[\|X^{q_n}_n(u)\|] + \E[\|X^{q_n}_n(u)\|^2]^{1/2} \Big)\diff u\\
		\label{eq:tight.X.1}
		&\quad + \|\sigma\|_\infty\bigg(\E\bigg[\int_s^t\|q_n(u)\|\diff u\bigg] + \frac{1}{\sqrt{n}}(t -s) \bigg).
	\end{align}
	Notice that the sequence $\E[\int_0^1\|q_n(t)\|\diff t]$ is bounded.
	In fact, by the coersivity condition on $g_1$, for each $c_1>0$, there is $N, c_2>0$ and $\varepsilon>1$ such that $g_1(t,q)\ge c_1\|q\|^\varepsilon$ for all $\|q\|\ge N$ and $g_1(t,q)\ge -c_2$ for all $(t,q)$.
	Thus, we always have 
	\begin{align}
	\notag
		\|q_n(t)\|^\varepsilon &\le \|q_n(t)\|^\varepsilon1_{\{\|q_n\|\ge N \}} + \|q_n(t)\|^\varepsilon1_{\{\|q_n\|< N \}}\\
			\label{eq:bound.q.coers}
		&\le  \frac{1}{c_1}g_1(t,  q_n(t) ) + N^\varepsilon \quad \text{for all } n\ge1,
	\end{align}
 	which shows by \eqref{eq:bound.g.lower} that $\E[\int_0^1\|q_n(t)\|\diff t]$ is bounded.
	In turn, this bound allows to conclude that $\E[\|X^{q_n}_n(u)\|^2]$ is bounded (this follows exactly as in the proof of \eqref{eq:bounded.moment} with $\|\varphi\|_1$ therein replaced by $\E[\int_0^1\|q_n(t)\|\diff t]$.
	Thus, by \eqref{eq:bounded.moment}, the first term on the right hand side of \eqref{eq:tight.X.1} is bounded by $C(t-s)$ for some $C>0$, so that due to \eqref{eq:bound.q.coers} we obtain
	\begin{align*}
		\E\Big[\|X^{q_n}_n(t) - X^{q_n}_n(s)\|\Big] &\le \Big(C + \frac{\|\sigma\|_\infty}{\sqrt{n}} \Big)(t-s) + \|\sigma\|_\infty\E\bigg[\int_s^t\frac{1}{ c_1}g_1(u,q_n(u)) + N^\varepsilon \diff u\bigg]^{1/\varepsilon},
	\end{align*}
	Thus, by \eqref{eq:bound.g.lower}, for any sequence $\delta_n\downarrow0$,
	\begin{align*}
		\lim_{n\to \infty}\sup_{n\ge1}\sup_{\tau}\E\Big[\|X^{q_n}_n(\tau + \delta_n) - X^{q_n}_n(\tau)\|\Big] \le \limsup_{n\to \infty} \Big(C + \frac{\|\sigma\|_\infty}{\sqrt{n}} \Big)\delta_n + \|\sigma\|_\infty\Big(\frac{C}{ c_1}\Big)^{1/\varepsilon},
	\end{align*}
	where the supremum on the left hand side is over stopping times $\tau$ with values in $[0,1 - \delta_n]$.
	Since $c_1>0$ was taken arbitrary, we conclude that
	\begin{align*}
		\lim_{n\to \infty}\sup_{n\ge1}\sup_{\tau}\E\Big[\|X^{q_n}_n(\tau + \delta_n) - X^{q_n}_n(\tau)\|\Big] =0.
	\end{align*} 
	Therefore, $(X^{q_n}_n)_{n\ge1}$ is tight in the set of continuous functions as a consequence of Aldous' tightness criterion \cite[Theorem 16.11]{kallenberg}.
	Thus, by Skorohod's representation theorem and passing to a subsequence, there is a probability space $(\Omega, \mathcal{F},\P)$ on which $(X^{q_n}_n)_{n\ge1}$ converges almost surely to some continuous process $X$, and that $(\mu^{q_n}_n(t))_{n\ge1}$ converges to $\mu(t)$ weakly.

	We now show that $X = \Phi^q$ (with $q$ the weak limit of $q_n$ introduced above).
	That is, we need to show that $X$ satisfies the equation
	\begin{equation}
	\label{eq:aux.eq.q}
	 	\diff X(t) = b\big( t, X(t), \mu(t) \big) + \sigma\big(t, X(t), \mu(t) \big)q(t)\diff t,\quad   
	 	X(s) = x.
	 \end{equation} 
	 First it is clear that $X(0) = x$.
	 By weak convergence of $(\mu^{q_n}_n(t))_{n\ge1}$, and the fact that $\EE\big[\|X^{q_n}_n(t)\|^2\big]$ converges to $\EE\big[\|X(t)\|^2\big]$, it follows that $(\mu^{q_n}_n(t))_{n\ge1}$ converges in second order Wasserstein distance, see \cite[Theorem 5.5]{MR3752669}. 
	 Hence, by continuity of $b$ and dominated convergence it holds that
	 \begin{equation}
	 \label{eq:conv.b.limit}
	 	\int_s^tb(u, X^{q_n}_n(u), \mu^{q_n}_n(u))\diff u \xrightarrow[n\to\infty]{}  \int_s^tb(u, X(u), \mu(u))\diff u\quad \P\text{--a.s.}
	 \end{equation}
	On the other hand, boundedness of $\sigma$ implies that $\frac{1}{\sqrt{n}}\int_0^t\sigma(u, X^{q_n}_n(u), \mu^{q_n}_n(u))\diff W(u)$ converges to zero $\P$-a.s. and by Lemma \ref{Lem:tightness},
	 \begin{equation}
	 \label{eq:conv.q.sigma}
	 	\int_s^t\sigma\big( u, X(u), \mu(u) \big) q_n(u)\diff u \xrightarrow[n\to\infty]{}  \int_s^t\sigma\big( u, X(u), \mu(u) \big) q(u)\diff u\quad \P\text{--a.s.}
	 \end{equation}
	Now, applying H\"older inequality, with $\widehat\varepsilon$ the H\"older conjugate of $\varepsilon$ in \eqref{eq:bound.q.coers}, we have
	\begin{align*}
		&\E\bigg[\int_s^t\Big(\sigma(u, X^{q_n}_n(u), \mu^{q_n}_n(u)) - \sigma(u, X(u), \mu(u)) \Big) q_n(u) \diff u\bigg] \\
		&\le \E\bigg[\int_0^1\|q_n(u)\|^\varepsilon\diff u\bigg]^{1/\varepsilon}\E\bigg[\int_0^1\|\sigma(u, X^{q_n}_n(u), \mu^{q_n}_n(u)) - \sigma(u, X(u), \mu(u))\|^{\widehat\varepsilon}\diff u\bigg]^{1/\widehat\varepsilon}\\
		&\le  \bigg( \frac{\|\sigma\|_\infty}{c_1}\E\bigg[\int_0^1 g_1(u,q_n(u))\diff u\bigg] + N\bigg)^{1/\varepsilon}\E\bigg[\int_0^1\|\sigma(u, X^{q_n}_n(u), \mu^{q_n}_n(u)) - \sigma(u, X(u), \mu(u))\|^{\widehat\varepsilon}\diff u\bigg]^{1/\widehat\varepsilon}\\
		&\le \bigg( \frac{\|\sigma\|_\infty}{c_1}C + N\bigg)^{1/\varepsilon}\E\bigg[\int_0^1\|\sigma(u, X^{q_n}_n(u), \mu^{q_n}_n(u)) - \sigma(u, X(u), \mu(u))\|^{\widehat\varepsilon}\diff u\bigg]^{1/\widehat\varepsilon},
	\end{align*}
	where the latter inequality follows by \eqref{eq:bound.g.lower}.
	This implies, using continuity and boundedness of $\sigma$ as well as dominated convergence, that
	\begin{equation}
	\label{eq:conv.sigma.zero}
		\E\bigg[\int_s^t\Big(\sigma(u, X^{q_n}_n(u), \mu^{q_n}_n(u)) - \sigma(u, X(u), \mu(u)) \Big) q_n(u) \diff u \bigg]  \xrightarrow[n\to\infty]{} 0 .
	\end{equation}
	Combine this with \eqref{eq:conv.q.sigma} to conclude that we can find a further subsequence such that
	\begin{equation}
	\label{eq:conv.sigma.limit}
		\int_s^t\sigma(u, X^{q_n}_n(u), \mu^{q_n}_n(u))q_n(u)\diff u  \xrightarrow[n\to\infty]{}   \int_s^t \sigma(u, X(u), \mu(u))q(u)\diff u \quad \P\text{-a.s}.
	\end{equation}
	Due to \eqref{eq:conv.b.limit} and continuity of $X$ this shows that $X$ satisfies \eqref{eq:aux.eq.q}, i.e. $X = \Phi^q$.

	Coming back to \eqref{eq:upper.bound.n.close}, we conclude by dominated convergence and \eqref{eq:fatou.upper} that
	\begin{align*}
		\limsup_{n\to \infty}\rho_s^{g_n}\big(F(X_n(1),\mu_n(1)) \big)
		&\le \E\bigg[ F(\Phi^q(1),\mu^q(1)) - \int_s^1g\big( t, q(t),\Phi^q(t),\mu^q(t) \big)\diff t  \bigg]\\
		&\le \sup_{\varphi\in \mathcal{H}}\bigg( F(\Phi(1),\delta_{\Phi(1)} ) - \int_s^1 g(t, \varphi(t),\Phi(t),\delta_{\Phi(t)})\diff t  \bigg),
	\end{align*}
	where we used the fact that $\mu(t) = \cL(\Phi(t)) = \delta_{\Phi(t)}$ when $\Phi(t)$ is deterministic.
	This concludes the proof.
\end{proof}
As a direct consequence of the above result, we derive the standard large deviation result in its Laplace principle form stated in the introduction.
\begin{proof}[Proof of Corollary \ref{cor:F-WThm}]
	The corollary follows from Theorem \ref{thm:abstract.schilder}, with the choice of function $g(t,q,x,\mu): = \frac12\|q\|^2$.
	In this case, the function $g_n$ introduced in the statement of Theorem \ref{thm:abstract.schilder} becomes $g_n(t,q,x,\mu) = \frac{1}{2n}\|q\|^2$ and we have
	\begin{equation*}
		\rho^{g_n}_0\big(F(X_n(1), \mu_n(1)) \big) = \frac1n\log \E\Big[e^{n F(X_n(1), \mu_n(1))} \Big].
	\end{equation*}
	Therefore, it holds that
	\begin{align*}
		\lim_{n\to \infty}\frac1n\log \E\Big[e^{n F(X_n(1), \mu_n(1))} \Big] &= \sup_{\varphi \in \mathcal{H}}\bigg(F\big(\Phi^\varphi(1), \delta_{\Phi^\varphi(1)}\big) - \frac12\int_0^1\|\varphi(t)\|^2\diff t \bigg)\\
		& = \sup_{\Phi \in \mathcal{H}_m}\Big(F\big(\Phi^\varphi(1), \delta_{\Phi^\varphi(1)}\big) - I(\Phi) \Big)
	\end{align*}
	with $I$ defined in the statement of the corollary.
\end{proof}
Let us conclude this subsection with a slight extension of Theorem \ref{thm:abstract.schilder}, namely to the random initial position case.
This extension will play a crucial role in the analysis of vanishing viscosity done in the ensuing subsection.
\begin{corollary}
\label{cor:F.W-random}
	Assume that the conditions $(A1)$, $(A2)$ and $(A3)$ are satisfied, and let $g_n$ be defined as in  Theorem \ref{thm:abstract.schilder}.
	Let $s \in [0,1]$ and $\xi \in \L^2(\R^\xdim, \cF_s)$ be given, and let $X_n$ solve the SDE
	\begin{equation}
	\label{eq:McKV.SDE.n.random}
		\begin{cases}
		dX_n(t) = b\big(t, X_n(t), \mu_n(t)\big)\diff t + \frac{1}{\sqrt{n}}\sigma\diff W(t),\\
		 \quad X_n(s) = \xi,\,\, \mu_n(t) = \text{law}(X_n(t)).
		\end{cases}
	\end{equation}
	For every $F \in C_b(\mathbb{R}^m\times \cP_2(\mathbb{R}^m))$ it holds
	\begin{equation*}
	\label{eq:abstract.schilder.determ}
		\lim_{n\to \infty}\rho^{g_n}_s\big(F(X_n(1), \mu_n(1))\big) = \esssup_{ q \in  \mathcal{L}}\E\bigg[F\big( \Phi^q(1),\mu^q(1) \big) - \int_s^1g\big(t, q(t), \Phi^q(t),\mu^q(t) \big)\diff t \bigg| \cF_s\bigg]
	\end{equation*}
	where $\Phi^q$ is the solution (when it exists) of the ordinary differential equation
	\begin{equation*}
	\label{eq:ODE.dertm}
		\begin{cases}
			\diff \Phi^\varphi(t)  =b\big(t, \Phi^q(t), \cL(\Phi^q(t)) \big) + \sigma\big(t, \Phi^q(t),\cL(\Phi^q(t)) \big)\varphi(t)\diff t\\
			\Phi^q(s) = \xi, \quad\mu^q(t) = \mathrm{law}(\Phi^q(t)).
		\end{cases}
	\end{equation*}
\end{corollary}
\begin{proof}
	The proof is essentially the same as that of Theorem \ref{thm:abstract.schilder}.
	In fact, the proof of the lower bound, i.e. the inequality "$\ge$" is unchanged, starting with at arbitrary $q\in \cL$ instead of $\varphi \in \mathcal{H}$ with $\|x\|^2$ in \eqref{eq:bounded.moment} replaced by $\EE[\|\xi\|^2]$.
	Notice that taking $\sigma$ constant allows to still get the bound \eqref{eq:boud.varphi.detem} in the present random control case.

	The proof of the upper bound differs only from the fact that we do not have access to a $1/n$-optimal control $(q_n)$ as in \eqref{eq:upper.bound.n.close}.
	To construct such a sequence, we need to first show that the set 
	$$\bigg\{\E\bigg[ F\big(X^{q}(1), \mu^{q}(1))\big) - \int_s^1g\big(u,q_n(u),X^{q_n}_n(u),\mu^{q_n}_n(u)\big)\diff u \bigg| \cF_s \bigg], q \in \cL \bigg\}$$
	is directed upward in the sense that for every $q_1,q_2 \in \cL$ there is $\hat q \in \cL$ such that $M_s^{\hat q} \ge \max\{M^{q_1}_s, M^{q_2}_s\}$, with
	\begin{equation*}
		M^{q}_s := \E\bigg[ F\big(X^{q}(1), \mu^{q}(1))\big) - \int_s^1g\big(u,q(u),X^{q}_n(u),\mu^{q}_n(u)\big)\diff u \bigg| \cF_s \bigg].
	\end{equation*}
	In fact, given $q_1, q_2\in \cL$, put $\tau := \inf\{t>s : M^{q_1}_t< M^{q_2}_t\}\wedge T$ and $\hat q:= q_11_{[0,\tau]}+ q_21_{[\tau, 1]}$.
	Then, it holds $M^{\hat q}_s \ge \max \{M^{q_1}_s, M^{q_2}_s\}$. 
	Thus, it follows by \cite[Theorem A.37]{FS3dr} that there is a sequence $(q_k)_{k\ge1}$ such that we have the increasing limit
	\begin{align}
		  \rho_s^{g_n}\big(F(X_n(1),\mu_n(1)) \big) = \lim_{k\to \infty}\E\bigg[ F\big(X^{q_k}_n(1), \mu^{q_k}_n(1))\big) - \int_s^1g\big(u,q_k(u),X^{q_k}_n(u),\mu^{q_k}_n(u)\big)\diff u \bigg],
	\end{align}
	thus by boundedness of $F$ and $g_2$ we have
	\begin{equation*}
		\E\bigg[\int_s^1g_1(u, q_k(u))\diff u \bigg| \cF_s\bigg] \le C.
	\end{equation*}
	Thus, using Lemma \ref{Lem:tightness} and arguing as in the proof of Theorem \ref{thm:abstract.schilder}, $(q_k)_{k\ge1}$ admits a subsequence converging to a process $q_n \in \cL$ and it holds
	\begin{equation*}
		\rho_s^{g_n}\big(F(X_n(1),\mu_n(1)) \big) = \E\bigg[ F\big(X^{q_n}_n(1), \mu^{q_n}_n(1))\big) - \int_s^1g\big(u,q_n(u),X^{q_n}_n(u),\mu^{q_n}_n(u)\big)\diff u \bigg| \cF_s \bigg].
	\end{equation*}
	This puts us exactly in the position of \eqref{eq:upper.bound.n.close}.
	The rest of the proof is the same.
\end{proof}

\subsection{Vanishing viscosity on the Wasserstein space}
\label{sec:PDE.proof}
In this subsection, we assume that $\sigma$ is constant and satisfies \eqref{eq:sigma.elliptic}.
On the way to drawing the link between the scaling limit theorems derived above and PDEs on the Wasserstein space, we will derive yet another formulation of the representation given in Theorem \ref{thm:var.rep.McKV}.
In the stochastic control language, we will show that the control problem on the right hand side of \eqref{eq:var.rep.abstract.McKV} has the same value when the set of admissible controls $\cL$ (which are open-loop) is replaced by the so-called Markovian controls.
Recall that $q$ is a Markovian control if there is a function $\varphi:[0,1]\times \R^\xdim\to \R^d$ such that $q(t) = \varphi(t,X^q_t)$.
That is, $q(t)$ is a function of the time-$t$ value of the state process.
The main argument used here to switch to Markovian controls will be the Mimicking theorem for McKean-Vlasov equations \cite{Lac-Shk-Zha2020}, see also \cite{Ben-Cam-DiPer,LackerSPA15} similar applications.
This view point of the control problem will allow to re-write it as a control of Fokker--Plank equations and thus to finally write the stochastic control problem as a deterministic one.
\begin{proposition}
\label{prop:abstract.schilder.deterministic}
	Assume that the conditions $(A1)$, $(A2)$ and $(A3)$ are satisfied.
	Let $X_n$ solve the SDE \eqref{eq:McKV.SDE.n.random} and let $g_n$ be defined by \eqref{eq:def.g.n}.
	Then for every $F \in C_b(\R^\xdim,\cP_2(\mathbb{R}^m))$ it holds
	\begin{equation}
	\label{eq:Abstract.Shilder.Deterministic}
		\lim_{n\to \infty}\rho^{g_n}_s(F(X_n(1),\mu_n(1)) = \sup_{\phi }\bigg(\widetilde F_s(\mu^\phi(1)) - \int_t^1\widetilde g_s\Big(u,\phi\big(u,\cdot, \mu^\phi(u)\big), \mu^\phi(u) \Big)\diff u \bigg)
	\end{equation}
	where the supremum is taken over Borel functions $\phi:[0,1]\times \R^\xdim \times \cP_2(\R^\xdim)\to \R^\bmdim$ such that $\phi(t, \cdot,\mu^\phi) \in \L^2(\R^\bmdim, \mu^\phi)$,
	with 
	$$\widetilde{F}_s(\mu) := \int_{\R^\xdim}F(	x,\mu)\mu(\diff x| \cF_s)\quad \text{and}\quad 
	\widetilde g_s(t, \phi(t,\cdot,\mu), \mu):=\int_{\R^\xdim}g(t, \phi(t,x,\mu),x,\mu)\mu(\diff x| \cF_s),$$
	where $\mu(\cdot|\cF_s)$ denotes the conditional distribution of $\mu$ given $\cF_s$
	and where $\mu^\phi$ satisfies the following continuity equation in the sense of distributions:
	\begin{equation}
	\label{eq:first.order.FP}
		\partial_t\mu^\phi = -\partial_x\Big\{\Big[ b(t, \cdot, \mu^\phi) + \sigma\phi(t, \cdot, \mu^\phi)\Big]\mu^\phi \Big\} .
	\end{equation}

	If $b$ does not depend on the law, that is, $b(t, x,\mu) = b(t,x)$, then it holds
		\begin{equation}
	\label{eq:Abstract.Shilder.Deterministic.b.neq.mu}
		\lim_{n\to \infty}\rho^{g_n}_s(F(X_n(1),\mu_n(1)) = \sup_{\phi}\bigg(\widetilde F_s(\mu^\phi(1)) - \int_t^1\widetilde g_s(u,\phi(u,\cdot), \mu^\phi(u))\diff u \bigg)
	\end{equation}
	where the supremum is taken over Borel functions $\phi:[0,1]\times \R^\xdim \to \R^\bmdim$ such that $\phi \in \L^2(\R^\bmdim, \mu^\phi)$,
	with $\mu^\phi$ satisfying the following continuity equation in the sense of distributions:
	\begin{equation}
	\label{eq:first.order.FP}
		\partial_t\mu^\phi = -\partial_x\Big\{\Big[ b(t, \cdot) + \sigma\phi(t, \cdot)\Big]\mu^\phi \Big\}.
	\end{equation}
\end{proposition}
\begin{proof}
	We showed in Corollary \ref{cor:F.W-random} that with probability one, 
	\begin{equation*}
		\lim_{n\to \infty}\rho^{g_n}_s\big(F(X_n(1), \mu_n(1))\big) = \esssup_{ q \in  \mathcal{L}}E\bigg[F(\Phi^q(1),\text{law}(\Phi^q(1))) - \int_s^1g\Big(t, q(t), \Phi^q(t),\text{law}(\Phi^q(t)) \Big)\diff t \bigg| \cF_s \bigg]
	\end{equation*}	
	where $\Phi^q$ solves the (random) SDE
	\begin{equation*}
		\Phi^q(t) = \xi + \int_s^tb\big(u, \Phi^q(u), \text{law}(\Phi^q(u)) \big) + \sigma q(u)\diff u.
	\end{equation*}
	Let $q$ be of the form 
	\begin{equation*}
		q(t) = \phi\big(t, \Phi^q(t), \text{law}(\Phi(t)) \big)
	\end{equation*}
	for some Borel--measurable function $\phi:[0,1]\times \R^\xdim \times \cP_2(\R^\xdim) \to \R^d$ that is Lipschitz--continuous in its last two arguments.
	Let us denote by $\mu^\phi(t)$ the law of $\Phi^q(t)$ with this specification of $q$.
	Fix a test function $f\in C^2_c(\R^\xdim)$, i.e. a twice continuously differentiable function with compact support.
	Applying It\^o's formula, we have
	\begin{align*}
		\int_s^1\langle f, \partial_t\mu^\phi(u)\rangle \diff u& = \E\Big[f(\Phi^q(1)) - f(\xi) \Big]
	 	 = \langle f, \mu^\phi(1)\rangle - \langle f, \mu^\phi(s)\rangle\\
	 	& = - \int_s^1\Big\langle f, \partial_x\Big\{\Big[ b(u, \cdot, \mu^\phi(u))+ \sigma\phi(u, \cdot, \mu^\phi(u))\Big] \Big\}\Big\rangle \diff u.
	 \end{align*}
	 This shows that $\mu^\phi$ satisfies \eqref{eq:first.order.FP}.
	 Therefore, we have 
	\begin{equation*}
		\lim_{n\to \infty}\rho^{g_n}_s(F(X_n(1),\mu_n(1)) \ge \sup_{\phi }\bigg(\widetilde F_s(\mu^\phi(1)) - \int_s^1\widetilde g_s(u,\phi(u,\cdot,\mu^\phi(u)), \mu^\phi(u))\diff u \bigg).
	\end{equation*}

	\vspace{.3cm}
	On the other hand, for every $\varepsilon>0$ as in the proof of Corollary \ref{cor:F.W-random}, 
	there is $q_\varepsilon \in \cL$ such that
	\begin{equation*}
		\lim_{n\to \infty}\rho^{g_n}_s\big(F(X_n(1), \mu_n(1))\big) \le \E\bigg[F\big( \Phi^{q_\varepsilon}(1), \text{law}(\Phi^{q_\varepsilon}(1) \big) - \int_s^1g\big(u,q_\varepsilon(u), \Phi^{q_\varepsilon}(u), \text{law}(\Phi^{q_\varepsilon}(u)) \big)\diff u \bigg|\cF_s \bigg] + \varepsilon.
	\end{equation*}
	By the mimicking theorem for McKean-Vlasov equations, see e.g. \cite[Corollary 1.6]{Lac-Shk-Zha2020} there is a probability space $(\widehat\Omega, \widehat\cF,\widehat\P)$ carrying a $d$-dimensional Brownian motion $W$ such that, putting 
	\begin{equation*}
		 \phi^\varepsilon(t, x,\mu) := \E\Big[q_\varepsilon(t)\Big| \Phi^{q_\varepsilon}(t) = x,\text{law}(\Phi^{q_\varepsilon}(t)) = \mu \Big],
	\end{equation*}
	there is $\widehat \Phi^\phi$ satisfying
	\begin{equation}
	\label{eq:equal.laws}
		\widehat \Phi^{\phi^\varepsilon}(t) \overset{ d} {= } \Phi^{q_\varepsilon}(t)
	\end{equation}
	and such that
	\begin{equation*}
		\diff \widehat \Phi^{\phi^\varepsilon}(t) = b\big(t, \widehat \Phi^{\phi^\varepsilon}(t),\widehat\mu^{\phi^\varepsilon}(t) \big) + \sigma\phi^\varepsilon\big(t, \widehat \Phi^{\phi^\varepsilon}(t), \widehat\mu^{\phi^\varepsilon}(t)\big)  \diff t 
	\end{equation*}
	with $\widehat \mu^{\phi^\varepsilon}(t)$ the law of $\widehat \Phi^{\phi^\varepsilon}(t)$ under $\widehat\P$.
	Therefore, by \eqref{eq:equal.laws} and convexity of $g$ in its second component and Jensen's inequality, we have
	\begin{align*}
		&\lim_{n\to \infty}\rho^{g_n}_s\big(F(X_n(1), \mu_n(1))\big)\\
		&\qquad \le \E^{\widehat\P}\bigg[F\big( \widehat \Phi^{\phi^\varepsilon}(1), \widehat\mu^{\phi^\varepsilon}(1) \big) - \int_s^1g\Big(u,\phi^\varepsilon\big(u, \widehat \Phi^{\phi^\varepsilon}(u), \widehat\mu^{\phi^\varepsilon}(u)\big), \widehat \Phi^{\phi^\varepsilon}(u), \widehat\mu^{\phi^\varepsilon}(u) \Big)\diff u \bigg|\cF_s \bigg] + \varepsilon\\
		&\qquad = \esssup_{\phi}\bigg(\widetilde F_s(\widehat\mu^\phi(1)) - \int_s^1\widetilde g_s\Big(u,\phi(u,\cdot, \widehat\mu^\phi(u)), \widehat\mu^\phi(u) \Big)\diff u \bigg) + \varepsilon
	\end{align*}
	with the supremum taken over Borel and integrable functions $\phi:[0,1]\times \R^\xdim \times \cP_2(\R^\xdim)\to \R^\bmdim$.
	This shows \eqref{eq:Abstract.Shilder.Deterministic} since $\varepsilon>0$ was chosen arbitrary.

\vspace{.3cm}

	The proof is exactly the same when $b$ does not depend on $\mu$, except that in the mimicking argument, it is enough to put
	\begin{equation*}
		 \phi^\varepsilon(t, x) := \E\Big[q_\varepsilon(t)\Big| \Phi^{q_\varepsilon}(t) = x \Big].
	\end{equation*}
\end{proof}

\begin{proof}[Proof of Theorem \ref{thm:PDE.convergence}]
	The main argument of the proof is already given in the proof of Proposition \ref{prop:abstract.schilder.deterministic}.
	Here, it remains only to show that the left hand side of Equation \eqref{eq:Abstract.Shilder.Deterministic} is in fact the (limit of the) solution of the PDE \eqref{eq:Wassers.PDE.n}. 
	Let $X_n$ solve the SDE \eqref{eq:McKV.SDE.n} and apply the extension of It\^o's formula to functions of laws of diffusions, see \cite[Theorem 5.104]{MR3752669}, we have
	\begin{align*}
		&\diff \cV_n(t, X_n(t), \mu_n(t))\\
		& = \partial_t\cV_n + \partial_x\cV_n\cdot b + \frac{1}{2n}\mathrm{Tr}\big[\partial_{xx}\cV_n\sigma\sigma^\top\big] +\int_{\mathbb{R}^\xdim}\partial_\mu \cV_n\big(t, X_n(t), \mu_n(t) \big)(a)b(t, a, \mu_n(t))\diff \mu_n(a)\\
		&  + \frac{1}{2n}\int_{\mathbb{R}^\xdim}\mathrm{Tr}\big[ \partial_a\partial_\mu \cV_n\big(t, X_n(t), \mu_n(t)\big)(a)\sigma\sigma^\top(t,a,\mu)\big]\diff \mu(a)\diff t + \frac{1}{\sqrt{n}}(\partial_x\cV_n\sigma) \big(t, X_n(t), \mu_n(t)\big)\cdot \diff W(t)
	\end{align*}	
	where, when the argument is not given, the function is evaluated at $(t, X_n(t), \mu_n(t))$.
	Since $\cV_n$ is a classical solution of the equation \eqref{eq:Wassers.PDE.n}, it thus follows that 
	$$
		(Y_n, Z_n) := \Big( \cV_n(t, X_n(t), \mu_n(t)), \frac{1}{\sqrt{n}} \partial_x\cV_n\sigma(t, X_n(t), \mu_n(t)) \Big)
	$$ solves the BSDE
	\begin{equation*}
		\begin{cases}
			\diff Y_n(t) = -f\big(t, \sqrt{n}Z_n(t), X_n(t), \mu_n(t) \big)\diff t + Z_n(t)\cdot\diff W(t)\\
			Y_n(1) = \cV(1, X_n(1), \mu_n(1)) = F(X_n(1), \mu_n(1)).
		\end{cases}
	\end{equation*}
	Thus, since $q\mapsto g_n(\cdot, q, \cdot, \cdot):= g(\cdot, q/\sqrt{n}, \cdot, \cdot)$ is the convex conjugate of the function $z\mapsto f(\cdot, \sqrt{n}z, \cdot,\cdot)$, it follows by \cite[Theorem 3.4]{tarpodual} that
	\begin{equation*}
		Y_n(t) = \rho^{g_n}_t\big( F(X_n(1), \mu_n(1)) \big),
	\end{equation*}
	showing that $\cV_n(t, X_n(t), \mu_n(t) ) = \rho^{g_n}_t\big( F(X_n(1), \mu_n(1)) \big)$.
	Hence, by Proposition \ref{prop:abstract.schilder.deterministic}, it follows that
	\begin{equation}
		\lim_{n\to \infty}\cV_n(s, \xi, \nu ) = \esssup_{\phi }\bigg(\widetilde F_s(\mu^\phi(1)) - \int_t^1\widetilde g_s\Big(u,\phi(u,\cdot, \mu^\phi(u)), \mu^\phi(u))\diff u \bigg).
	\end{equation}
	with the supremum taken over integrable functions $\phi:[0,1]\times \R^\xdim\times \cP_2(\R^\xdim) \to \R^\bmdim$, and when $b$ does not depend on $\mu$, the supremum is taken over integrable functions $\phi:[0,1]\times \R^\xdim\to \R^\bmdim$ as argued in Proposition \ref{prop:abstract.schilder.deterministic}.
	This concludes the proof.
\end{proof}

\section{A noteworthy consequence}
In this final section, we present another application of our variational representation result.
It concerns an application to functional inequalities.

In this subsection, assume that $\sigma$ satisfies $(A3)$, i.e. it is constant and satisfies the uniform ellipticity condition \eqref{eq:sigma.elliptic}, and that the function $b$ is linear in $x$ and $\mu$.
That is, there are bounded measurable functions $\alpha,\beta,\gamma:[0,T]\to \R$ such that
\begin{equation*}
	b(t, x, \mu) = \alpha(t) + \beta(t) x + \gamma(t)\int_{\R^\xdim}x\mu(\diff x).
\end{equation*}
Under these conditions, the SDE
\begin{equation}
\label{eq:SDE.prekopa.leindler}
	\diff X(t) = b(t, X_t, \mu(t))\diff t + \sigma\diff W(t),\quad X(0) = x, \quad \mu(t) = \text{law}(X_t)
\end{equation}
is well-posed \cite{MR3752669}.
We will show below that the law $\mu(t)$ of $X(t)$ satisfies Pr\'ekopa--Leindler inequality \eqref{eq:P.L.inequ}.
This is an integral version of the celebrated Brunn--Minkowski inequality \cite{Gardner02}, and it directly implies that the law $\mu(t)$ is log-concave, which is a key property for instance to allow efficient sampling algorithms using Langenvin monte--carlo methods.
\begin{proposition}[Pr\'ekopa-Leindler inequality]
\label{eq:prekopa.leindler}
	Let $t \in [0,1]$ and denote by $\mu_t$ the law of $X(t)$ in \eqref{eq:SDE.prekopa.leindler}.
	Assume that $g$ satisfies $(A2)'$ and depends only on $(t,q)$, i.e $g(t,q,x,\mu) = g(t,q)$.
	Let $0< \lambda < 1$ and $\ell_1,\ell_2$ and $\ell_3$ be three non-negative functions mapping $\R^\xdim$ to $\R$, belonging to $\L^1(\mu_t)$ and such that
	\begin{equation}
	\label{eq:P.L.condition}
	 	\ell_3((1-\lambda)x + \lambda y) \ge \ell_1(x)^{1-\lambda}\ell_2(y)^{\lambda},
	 \end{equation} 
	 for all $x, y \in \R^\xdim$.
	 Then,
	 \begin{equation}
	 \label{eq:p.l.abstract}
		\rho^g_0\big( \ell_3(X(t)) \big) \ge (1-\lambda)\rho^g_0\big( \ell_1(X(t)) \big) + \lambda \rho^g_0\big( \ell_2(X(t)) \big).
	 \end{equation}
	 In particular, we have Pr\'ekopa--Leindler inequality
	 \begin{equation}
	 \label{eq:P.L.inequ}
	 		\int_{\R^\xdim}\ell_3(x)\mu_t(\diff x) \ge \Big(\int_{\R^\xdim} \ell_1(x)\mu_t(\diff x) \Big)^{1-\lambda}\Big(\int_{\R^\xdim} \ell_2(x)\mu_t(\diff x) \Big)^{\lambda}.
	 \end{equation}
\end{proposition}
\begin{proof}
	The proof follows the standard stochastic approach of \citet{Bo00}.
	In fact, specializing the representation of Theorem \ref{thm:var.rep.McKV} to the case $g(t, q, x, \mu) = g(t,q)$ and $F(x, \mu) = F(x)$, we have
	\begin{equation*}
		\rho^g_0\big(F(X(t)) \big) = \sup_{q \in \cL}\E\bigg[F(X^q(t)) -\int_0^1g(u,q(u))\diff u\bigg].
	\end{equation*}
	Therefore for any $\varepsilon>0$, applying this inequality to $F(x):= \ell_1(x)$ and $F(x) = \ell_2(x)$, there are $q_1$ and $q_2$ such that
	\begin{align*}
		&(1-\lambda)\rho^g_0\big(\ell_1(X(t)) \big) + \lambda \rho^g_0\big(\ell_2(X(t)) \big)\\
		&\quad \le (1-\lambda)\E\bigg[\ell_1(X^{q_1}(t))  -\int_0^1g(u,q_1(u)) \diff u   \bigg] + \lambda\E\bigg[\ell_2(X^{q_1}(t)) - \int_0^1g(u,q_2(u))\diff u   \bigg] + \varepsilon\\
		&\quad \le \E\bigg[\ell_3\big((1-\lambda) X^{q_1}(t) + \lambda X^{q_2}(t) \big) - \int_0^1g(u,(1-\lambda)q_1(u) + \lambda q_2(u))\diff u   \bigg] + \varepsilon,
	\end{align*}
	where we used Jensen's inequality and \eqref{eq:P.L.condition}.
	Since $(1-\lambda) X^{q_1}+\lambda X^{q_2} = X^{(1-\lambda) q_1 + \lambda q_2 }$, we then have
	\begin{align*}
		&(1-\lambda)\rho^g_0\big(\ell_1(X(t)) \big) + \lambda \rho^g_0\big(\ell_2(X(t)) \big)\\
		&\le \E\bigg[\ell_3\big( X^{(1-\lambda)q_1 + \lambda q_2}(t) \big) - \int_0^1 g\big( u,(1-\lambda)q_1(u) + \lambda q_2(u) \big)\diff u   \bigg] + \varepsilon\\
		&\le \sup_{q \in \cL}\E\bigg[\ell_3(X^q(t))) - \int_0^1g(u,q(u))\diff u \bigg] + \varepsilon  = \rho^g_0\big( \ell_3(X(t))\big) + \varepsilon.
	\end{align*}
	Since $\varepsilon$ was taken arbitrary, this yields \eqref{eq:p.l.abstract}.
	This inequality allows to obtain \eqref{eq:P.L.inequ} by taking $g(t,q):=\frac12\|q\|^2$ since in this case $\rho^{g}_0(F(X(t))) = \log\E[e^{F(X(t))}]$ so that \eqref{eq:p.l.abstract} becomes
	\begin{equation*}
		\log\Big(\E[e^{\ell_3(X(t))}] \Big) \ge (1-\lambda)\log\Big(\E[e^{\ell_1(X(t))}] \Big) + \lambda \log\Big(\E[e^{\ell_2(X(t))}] \Big).
	\end{equation*}
	Taking exponential on both sides leads to \eqref{eq:P.L.inequ}.	
\end{proof}

\bibliographystyle{abbrvnat}
%\bibliography{references-Funct-ineq-conv}

\begin{thebibliography}{53}
\providecommand{\natexlab}[1]{#1}
\providecommand{\url}[1]{\texttt{#1}}
\expandafter\ifx\csname urlstyle\endcsname\relax
  \providecommand{\doi}[1]{doi: #1}\else
  \providecommand{\doi}{doi: \begingroup \urlstyle{rm}\Url}\fi

\bibitem[Ambrosio and Feng(2014)]{Ambro-Feng}
L.~Ambrosio and J.~Feng.
\newblock On a class of first order {H}amilton--{J}acobi equations in metric
  spaces.
\newblock \emph{J. Differential Equations}, 256\penalty0 (7):\penalty0
  2194--2245, 2014.

\bibitem[Backhoff-Veraguas et~al.(2020)Backhoff-Veraguas, Lacker, and
  Tangpi]{BaLaTa}
J.~Backhoff-Veraguas, D.~Lacker, and L.~Tangpi.
\newblock Nonexponential {S}anov and {S}childer theorems on {W}iener space:
  {BSDEs}, {S}chr\"odinger problems and control.
\newblock \emph{Ann. Appl. Probab.}, 30\penalty0 (3):\penalty0 1321--1367,
  2020.

\bibitem[Benazzoli et~al.(2020)Benazzoli, Campi, and {Di
  Persio}]{Ben-Cam-DiPer}
C.~Benazzoli, L.~Campi, and L.~{Di Persio}.
\newblock Mean field games with controlled jump--diffusion dynamics:
  {E}xistence results and an illiquid interbank marrket model.
\newblock \emph{Stoch. Proc. Appl.}, 130\penalty0 (11):\penalty0 6927--6964,
  2020.

\bibitem[Borell(2000)]{Bo00}
C.~Borell.
\newblock Diffusion equations and geometric inequalities.
\newblock \emph{Potential Anal.}, 12\penalty0 (1):\penalty0 49--71, 2000.

\bibitem[Bou\'e and Dupuis(1998)]{Boue-Dup}
M.~Bou\'e and P.~Dupuis.
\newblock A variational representation for certain functionals of {B}rownian
  motion.
\newblock \emph{Ann. Probab.}, 26\penalty0 (4):\penalty0 1641--1659, 1998.

\bibitem[{Capuzzo Dolcetta} and Leoni(2000)]{Cap-Leoni}
I.~{Capuzzo Dolcetta} and F.~Leoni.
\newblock On the vanishing viscosity approximation of a time dependent
  {H}amilton--{J}acobi equation.
\newblock In \emph{Recent Trends in Nonlinear Analysis}, pages 59--75.
  Birkh\"auser, 2000.

\bibitem[Cardaliaguet et~al.(2019)Cardaliaguet, Delarue, Lasry, and
  Lions]{carda15}
P.~Cardaliaguet, F.~Delarue, J.-M. Lasry, and P.-L. Lions.
\newblock \emph{The Master Equation and the Convergence Problem in Mean-Field
  Game}.
\newblock Princeton University Press, 2019.

\bibitem[Carmona and Delarue(2018{\natexlab{a}})]{MR3752669}
R.~Carmona and F.~Delarue.
\newblock \emph{Probabilistic theory of mean field games with applications.
  {I}}, volume~83 of \emph{Probability Theory and Stochastic Modelling}.
\newblock Springer, Cham, 2018{\natexlab{a}}.
\newblock Mean field FBSDEs, control, and games.

\bibitem[Carmona and Delarue(2018{\natexlab{b}})]{MR3753660}
R.~Carmona and F.~Delarue.
\newblock \emph{Probabilistic theory of mean field games with applications.
  {II}}, volume~84 of \emph{Probability Theory and Stochastic Modelling}.
\newblock Springer, Cham, 2018{\natexlab{b}}.
\newblock Mean field games with common noise and master equations.

\bibitem[Cecchin and Delarue(2020)]{Cecc-Del20}
A.~Cecchin and F.~Delarue.
\newblock Selection by vanishing common noise for potential finite state mean
  field games.
\newblock \emph{Preprint}, 2020.

\bibitem[Chassagneux et~al.(2014)Chassagneux, Crisan, and
  Delarue]{ChassagneuxCrisanDelarue_Master}
J.-F. Chassagneux, D.~Crisan, and F.~Delarue.
\newblock A probabilistic approach to classical solutions of the master
  equation for large population equilibria.
\newblock \emph{Forthcoming in Memoirs of the AMS}, 2014.

\bibitem[Crandall et~al.(2000)Crandall, Kocan, and \'Swiech]{Grandall-Koc-Sci}
M.~G. Crandall, M.~Kocan, and A.~\'Swiech.
\newblock $l^p$--{T}heory for fully nonlinear uniformly parabolic equations.
\newblock \emph{Comm. Part. Diff. Equa.}, 25\penalty0 (11-12):\penalty0
  1997--2053, 2000.

\bibitem[Delarue and Tchuendom(2020)]{Del-Tcheu20}
F.~Delarue and R.~F. Tchuendom.
\newblock Selection of equilibria in a linear quadratic mean-field game.
\newblock \emph{Stoch. Proc. Appl.}, 130\penalty0 (2):\penalty0 1000--1040,
  2020.

\bibitem[Delbaen et~al.(2011)Delbaen, Hu, and Bao]{Delbaen11}
F.~Delbaen, Y.~Hu, and X.~Bao.
\newblock Backward {SDEs} with superquadratic growth.
\newblock \emph{Probab. Theory Relat. Fields}, 150:\penalty0 145--192, 2011.

\bibitem[{dos Reis} et~al.(2019){dos Reis}, Salkeld, and
  Tugaut]{dReis-Sal-Tug18}
G.~{dos Reis}, W.~Salkeld, and J.~Tugaut.
\newblock Freidlin-{W}entzell {LDP} in path space for {McKean-Vlasov} equations
  and the functional iterated logarithm law.
\newblock \emph{Ann. Appl. Probab.}, 29\penalty0 (3):\penalty0 1487--1540,
  2019.

\bibitem[Drapeau et~al.(2013)Drapeau, Heyne, and Kupper]{DHK1101}
S.~Drapeau, G.~Heyne, and M.~Kupper.
\newblock Minimal supersolutions of convex {BSDEs}.
\newblock \emph{Ann. Probab.}, 41\penalty0 (6):\penalty0 3697--4427, 2013.

\bibitem[Drapeau et~al.(2016)Drapeau, Kupper, {Rosazza Gianin}, and
  Tangpi]{tarpodual}
S.~Drapeau, M.~Kupper, E.~{Rosazza Gianin}, and L.~Tangpi.
\newblock Dual representation of minimal supersolutions of convex {BSDEs}.
\newblock \emph{Ann. Inst. H. Poincar\'e Probab. Statist.}, 52\penalty0
  (2):\penalty0 868--887, 2016.

\bibitem[Dupuis and Ellis(2011)]{dupuis-ellis}
P.~Dupuis and R.~Ellis.
\newblock \emph{A weak convergence approach to the theory of large deviations},
  volume 902.
\newblock John Wiley \& Sons, 2011.

\bibitem[El~Karoui et~al.(1997)El~Karoui, Peng, and Quenez]{karoui01}
N.~El~Karoui, S.~Peng, and M.~C. Quenez.
\newblock Backward stochastic differential equations in finance.
\newblock \emph{Math. Finance}, 1\penalty0 (1):\penalty0 1--71, 1997.

\bibitem[Evans(1998)]{Evans1998}
L.~C. Evans.
\newblock \emph{Partial Differential Equations}, volume~19 of \emph{Graduate
  Studies in Mathematics}.
\newblock American Mathematical Society, Providence, RI, 1998.
\newblock ISBN 0-8218-0772-2.

\bibitem[Fleming(1966)]{Fleming_JMA}
W.~H. Fleming.
\newblock Duality and a priori estimates in {M}arkovian optimization problems.
\newblock \emph{J. Math. Anal. Appl.}, 16:\penalty0 254--279; Erratum, 19
  (1969), p. 204, 1966.

\bibitem[Fleming(1969)]{Fleming-JDE}
W.~H. Fleming.
\newblock The {C}auchy problem for a nonlinear first--order partial
  differential equation.
\newblock \emph{J. Differential Equation}, 5:\penalty0 515--530, 1969.

\bibitem[Fleming(1971)]{Fleming71}
W.~H. Fleming.
\newblock Stochastic control for small noise intensities.
\newblock \emph{SIAM J. Control Optim.}, 9\penalty0 (3):\penalty0 473--517,
  1971.

\bibitem[Fleming(1977/78)]{Fl78}
W.~H. Fleming.
\newblock Exit probabilities and optimal stochastic control.
\newblock \emph{Appl. Math. Optim.}, 4\penalty0 (4):\penalty0 329--346,
  1977/78.

\bibitem[Fleming and Soner(2006)]{Flem-Soner-second}
W.~H. Fleming and H.~M. Soner.
\newblock \emph{Controlled {M}arkov processes and viscosity solutions},
  volume~25 of \emph{Stochastic Modelling and Applied Probability}.
\newblock Springer, New York, second edition, 2006.

\bibitem[Fleming and Souganidis(1986)]{Fleming-Souga86}
W.~H. Fleming and P.~E. Souganidis.
\newblock Asymptotic series and the method of vanishing viscosity.
\newblock \emph{Indiana University Mathematics Journal}, 35\penalty0
  (2):\penalty0 425--447, 1986.

\bibitem[F{\"o}llmer and Schied(2011)]{FS3dr}
H.~F{\"o}llmer and A.~Schied.
\newblock \emph{Stochastic finance}.
\newblock Walter de Gruyter \& Co., Berlin, 3rd edition edition, 2011.
\newblock An introduction in discrete time.

\bibitem[Gangbo and \'Swiech(2015)]{Gang-Swe15}
W.~Gangbo and A.~\'Swiech.
\newblock Existence of a solution to an equation arising from the theory of
  mean field games.
\newblock \emph{J. Differential Equations}, 259\penalty0 (11):\penalty0
  6573--6643, 2015.

\bibitem[Gangbo and Tudorascu(2019)]{Gang-Tud19}
W.~Gangbo and A.~Tudorascu.
\newblock On differentiability in the {W}asserstein space and well-posedness
  for {H|}amilton--{J}acobi equations.
\newblock \emph{J. Math. Pures Appl.}, 125:\penalty0 119--174, 2019.

\bibitem[Gangbo et~al.(2008)Gangbo, Nguyen, and Tudorascu]{Gangbo08}
W.~Gangbo, T.~Nguyen, and A.~Tudorascu.
\newblock Hamilton--{J}acobi equations in the {W}assertein space.
\newblock \emph{Methods Appl. Anal.}, 15\penalty0 (2):\penalty0 155--183, 2008.

\bibitem[Gangbo et~al.(2009)Gangbo, Nguyen, and Tudorascu]{Gang-Ngu09}
W.~Gangbo, T.~Nguyen, and A.~Tudorascu.
\newblock Euler-{P}oisson systems as action--minimizing paths in the
  {W}asserstein space.
\newblock \emph{Arch. Ration. Mech. Anal.}, 192\penalty0 (3):\penalty0
  419--452, 2009.

\bibitem[Gardner(2002)]{Gardner02}
R.~Gardner.
\newblock The {B}runn-minkowski inequality.
\newblock \emph{Bull. Amer. Math. Soc.}, 39\penalty0 (9):\penalty0 355--405,
  2002.

\bibitem[Herrmann et~al.(2008)Herrmann, Imkeller, and
  Peithmann]{Herr-Imkeller-Pei08}
S.~Herrmann, P.~Imkeller, and D.~Peithmann.
\newblock Large deviation and a {K}ramers' typee law for self-stabilizing
  diffusions.
\newblock \emph{Ann. Appl. Probab.}, 18\penalty0 (4):\penalty0 1379--1423,
  2008.

\bibitem[Huang et~al.(2007)Huang, Caines, and Malham{\'e}]{Huang2007}
M.~Huang, P.~E. Caines, and R.~P. Malham{\'e}.
\newblock Large-population cost-coupled {LQG} problems with nonuniform agents:
  individual-mass behavior and decentralized {$\epsilon$}-{N}ash equilibria.
\newblock \emph{IEEE Trans. Automat. Control}, 52\penalty0 (9):\penalty0
  1560--1571, 2007.

\bibitem[Hynd and Kim(2015)]{HyndJFA}
R.~Hynd and H.~K. Kim.
\newblock Value functions in the wasserstein space: finite time horizons.
\newblock \emph{J. Funct. Anal.}, 269:\penalty0 968--997, 2015.

\bibitem[Kallenberg(2006)]{kallenberg}
O.~Kallenberg.
\newblock \emph{Foundations of modern probability}.
\newblock Springer Science \& Business Media, 2006.

\bibitem[Kobylanski(2000)]{kobylanski01}
M.~Kobylanski.
\newblock Backward stochastic differential equations and partial differential
  equations with quadratic growth.
\newblock \emph{Ann. Probab.}, 28\penalty0 (2):\penalty0 558--602, 2000.

\bibitem[Krylov(1980)]{KrylovBook}
N.~V. Krylov.
\newblock \emph{Controlled Diffusion Processes}.
\newblock Springer Verlag, 1980.

\bibitem[Lacker(2015)]{LackerSPA15}
D.~Lacker.
\newblock Mean field games via controlled martingale problems: {E}xistence of
  {M}arkovian equilibria.
\newblock \emph{Stoch. Proc. Appl.}, 125\penalty0 (7):\penalty0 2856--2894,
  2015.

\bibitem[Lacker et~al.(2020)Lacker, Shkolnikov, and Zhang]{Lac-Shk-Zha2020}
D.~Lacker, M.~Shkolnikov, and J.~Zhang.
\newblock Superposition and mimicking theorems for conditional {McKean-Vlasov}
  equations.
\newblock \emph{Preprint}, 2020.

\bibitem[Lasry and Lions(2007)]{MR2295621}
J.-M. Lasry and P.-L. Lions.
\newblock Mean field games.
\newblock \emph{Jpn. J. Math.}, 2\penalty0 (1):\penalty0 229--260, 2007.

\bibitem[Lehec(2013)]{Lehec}
J.~Lehec.
\newblock Representation formula for the entropy and functional inequalities.
\newblock \emph{Ann. Inst. H. Poincar\'e Probab. Statist.}, 49\penalty0
  (3):\penalty0 885--899, 2013.

\bibitem[Leindler(1972)]{Leindler72}
L.~Leindler.
\newblock On a certain converse of {H}\"older's inequality. {II}.
\newblock \emph{Acta Sci. Math. (Szeged)}, 33:\penalty0 217--223, 1972.

\bibitem[Lions(2007-2012)]{PLLcollege}
P.-L. Lions.
\newblock Cours du {C}oll{\`e}ge de {F}rance.
\newblock http://www.college-de-france.fr/default/EN/all/equ$_-$der/,
  2007-2012.

\bibitem[Pardoux and Peng(1992)]{Pardoux-Peng92}
E.~Pardoux and S.~Peng.
\newblock Backward stochastic differential equations and quasilinear parabolic
  partial differential equations.
\newblock In \emph{Stochastic partial differential equations and their
  applications}, volume 176 of \emph{Lecture Notes in Control and Inform. Sci.}
  Springer, Berlin, 1992.

\bibitem[Pr\'ekopa(1971)]{Prekopa71}
A.~Pr\'ekopa.
\newblock Logarithmic concave measures with application to stochastic
  programming.
\newblock \emph{Acta Sci. Math. (Szeged)}, 32:\penalty0 301--316, 1971.

\bibitem[Rockafellar(1976)]{rockafellar03}
R.~T. Rockafellar.
\newblock {Integral Functionals, Normal Integrands and Measurable Selections}.
\newblock In J.~Gossez, E.~Lami~Dozo, J.~Mawhin, and L.~Waelbroeck, editors,
  \emph{Nonlinear Operators and the Calculus of Variations}, volume 543 of
  \emph{Lecture Notes in Mathematics}, pages 157--207. Springer Berlin /
  Heidelberg, 1976.

\bibitem[Rockafellar and Wets(1998)]{rockafellar02}
R.~T. Rockafellar and R.~J.-B. Wets.
\newblock \emph{Variational Analysis}.
\newblock Springer, Berlin, New York, 1998.

\bibitem[Tugaut(2016)]{Tug16}
J.~Tugaut.
\newblock A simple proof of a {K}ramers' type law for self-stabilizing
  diffusions.
\newblock \emph{Electron. Commun. Probab.}, 21\penalty0 (11):\penalty0 1--7,
  2016.

\bibitem[van Handel(2018)]{vanHandel17}
R.~van Handel.
\newblock The {B}orell-{E}hrhard game.
\newblock \emph{Probab. Theory Related Fields}, 170:\penalty0 555--585, 2018.

\bibitem[Veretennikov(1981)]{Verete_USSR}
A.~Y. Veretennikov.
\newblock On strong solutions and explicit formulas for solutions of stochastic
  integral equations.
\newblock \emph{Math. USSR Sbornik}, 39\penalty0 (3), 1981.

\bibitem[Waagan(2008)]{Waagan}
K.~Waagan.
\newblock Convergence rate of monotone numerical scheme for
  {H}amilton--{J}acobi equations with weak boundary conditions.
\newblock \emph{SIAM J. Numerical Analysis}, 46\penalty0 (5):\penalty0
  2371--2392, 2008.

\bibitem[Zhang(2017)]{zhangbook}
J.~Zhang.
\newblock \emph{Backward Stochastic Differential Equations -- from linear to
  fully nonlinear theory}.
\newblock Springer, 2017.

\end{thebibliography}

\end{document}